\documentclass[12pt]{amsart}
\usepackage{amsfonts,amsmath,amssymb,comment,latexsym}
\numberwithin{equation}{section}

\newtheorem{Lem}{Lemma}[section]
\newtheorem{Prop}[Lem]{Proposition}
\newtheorem{Cor}[Lem]{Corollary}

\newtheorem{Thm}[Lem]{Theorem}

\theoremstyle{definition}
\newtheorem{Def}[Lem]{Definition}

\renewcommand{\k}{{\Bbbk}}
\renewcommand\o{\otimes}

\newcommand\D{{\mathcal D}}
\newcommand\Dg[1]{{\langle #1\rangle}}
\newcommand{\G}{{\widehat G}}
\newcommand{\GG}{{\widehat\Gamma}}
\newcommand{\ua}{\underline{a}}
\newcommand{\uc}{\underline{c}}
\newcommand{\ud}{\underline{d}}

\newcommand{\ui}{\underline {i}}
\newcommand{\uj}{\underline {j}}
\newcommand{\uo}{\underline {0}}
\renewcommand\ui{\underline{i}}
\renewcommand\uj{\underline{j}}
\newcommand\uk{\underline{k}}
\newcommand\ul{\underline{l}}
\newcommand\um{\underline{m}}
\newcommand{\un}{\underline {n}}

\newcommand\uN{\underline{N}}
\newcommand{\up}{\underline {p}}
\newcommand\ur{\underline{r}}
\newcommand\us{\underline{s}}
\newcommand\ut{\underline{t}}
\newcommand\uu{\underline{u}}
\newcommand\uv{\underline{v}}
\newcommand\ux{\underline{x}}
\newcommand\uy{\underline{y}}

\newcommand\lh{\rightharpoonup}
\newcommand\rh{\leftharpoonup}

\newcommand{\pf}{\medskip\noindent{\sc Proof: }}

\begin{document}

\title[Representations of Hopf Algebras]{Representations of Pointed Hopf Algebras and their Drinfel'd Quantum Doubles}
\author{Leonid Krop}
\thanks{Research by the first author partially supported by a grant from the College of Liberal Arts and Sciences at DePaul University}
\address{DePaul University\\Department of Mathematical Sciences\\2320 N. Kenmore\\Chicago, IL 60614}\email{lkrop@depaul.edu}
\author{David E. Radford}
\thanks{Research by the second author partially supported by NSA Grant H98230-04-1-0061}
\address {University of Illinois at Chicago\\Department of Mathematics, Statistics and Comuter Science\\801 South Morgan Street\\Chicago, IL 6o607-7045}\email{radford@uic.edu}

\begin{abstract}We study representations of nilpotent type nontrivial liftings of quantum linear spaces and their Drinfel'd quantum doubles. We construct a family of Verma- type modules in both cases and prove a parametrization theorem for simple modules. We compute the Loewy and socle series of Verma modules under a mild restriction on the datum of a lifting. We find bases and dimensions of simple modules. 
\end{abstract}
\date{\today}
\maketitle

\section*{Introduction}\label{intro}

Let $H$ be a finite-dimensional Hopf algebra over a field $\Bbbk$. The goal of this paper is two-fold. First, we want to describe structure of a family of Verma-type $H$-modules, when $H$ is a certain lifting of a quantum linear space, which entails determination of all simple $H$-modules. In the second place we carry out a similar program for the Drinfel'd quantum double of $H$.

We survey some previous related work. There has been significant interest in recent years in representation theory of nonsemisimple Hopf algebras and their quantum doubles \cite{{AB},{Ya},{Ch},{EGST},{KR1},{KR2},{CK},{RS}}. In the general setting of such algebra the primary focus is on classifying all simple $H$-modules in terms of simple modules of the coradical $H_0$ of $H$ \cite{{Ra2},{KR2}}. When $H_0=\Bbbk G$ where $G=G(H)$ is the group of grouplike elements of $H$, and $G$ is abelian with $\Bbbk$ a splitting field for $G$ of characteristic zero, the simple $H_0$-modules are given by the elements of the dual group $\G$. The problem of establishing a bijection between $\G$ and the set of isomorphism classes of simple $H$-modules will be called parametrization of simple $H$-modules. A typical example of pointed Hopf algebra with abelian group of grouplike is provided by the fundamental classification result of Andruskiewitsch and Schneider \cite{AS3}. Recently Radford and Schneider \cite{RS} proved the parametrization property for every algebra of the form $u({\mathcal D},\lambda)$ from the classification. This result generalizes an earlier theorem of Lusztig involving his small quantum groups \cite{Lu2}. The aforementioned theorems invariably involve a Hopf algebra $A$ with a triangular decomposition 
\begin{equation}\label{triangulatedalgebra} A=A^-\otimes A_0\otimes A^+\end{equation}
satisfying the following conditions. 
\begin{align}\label{normality} &A_0\, \text{is a subHopfalgebra},\, A^-\, \text{and}\, A^+ \text{are subalgebras}\\& \text{stable under}\, \text{ad}_{\ell} - \text{and}\, \text{ad}_r -\text{action of}\, A_0,\; \text{and}\nonumber\end{align}
 
\begin{equation}\label{nearlynilpotent} A^-=\Bbbk\oplus J^-,\; A^+=\Bbbk\oplus J^+\end{equation}
where $J^-,\,J^+$ are nilpotent ideals of $A^-,\,A^+$, respectively.

Let $H_i$ denote the $i$th term of the coradical filtration of $H$. By e.g. \cite{KR1} $H_1$ is a free $H_0$-module on a basis $\{1,x_1,\ldots,x_n\}$ where $x_i$ are some skew-primitive elements. When $H$ is generated by $H_1$ we say that $n$ is the {\em rank} of $H$. The simplest family of Hopf algebras within the class of pointed Hopf algebras $H$ with abelian $G(H)$ are liftings of quantum linear spaces completely described in \cite{AS1}. These Hopf algebras are natural generalizations of Lusztig's small quantum group $u_q(\mathfrak sl_2)$ associated to the simple Lie algebra of rank $1$. The concept of (Hopf) rank is applicable to them, and in small ranks $n=1,2$ the regular representation of $H$ is described in \cite{{KR1},{CK}}. For special cases of rank~$2$ liftings the entire finite-dimensional representation theory has been obtained in \cite{{Ch},{EGST}}.

We give an outline of the results. Let $H$ be a lifting of a quantum linear space, and $x_1,\ldots,x_n$ be skew-primitive generators of $H$. We say that $x_i$ and $x_j$ are {\em linked} if $x_ix_j-q_{ij}x_jx_i\neq 0$ for a certain root of unity $q_{ij}$. Let $\Gamma$ be the simple graph on the vertex set $\{1,2,\ldots,n\}$ with the edge set made up by pairs $(i,j)$ such that $x_i,x_j$ are linked. We say that $\Gamma$ is simply linked if every two vertices are connected by at most one edge. A lifting $H$ is called of {\em nilpotent type} if all $x_i$ are nilpotent.

The first part of the paper is concerned with algebras satisfying conditions \eqref{triangulatedalgebra}-\eqref{nearlynilpotent}, especially nilpotent simply linked lifting $H$ of a quantum linear space. In the latter case  $H_0=\Bbbk G(H)$. By analogy with the representation theory of $u_q(\mathfrak sl_2)$ we can use either subHopfalgebra $H^{\ge 0}:=H_0H^+$ or $H^{\le 0}:=H_0H^-$ to construct Verma-type modules $Z(\gamma)$ where $\gamma$ runs over $\G$. From the general Theorem \ref{generalparametrization} we have that $Z(\gamma)$ has a unique maximal submodule $R(\gamma)$, the radical of $Z(\gamma)$. Iterating this procedure gives the radical filtration $\{R^m(\gamma)\}$ of $Z(\gamma)$. The first major result of the paper is a description of $R^m(\gamma)$ and the (Loewy) layers of the filtration. We also show that the socle filtration of $Z(\gamma)$ coincides with the radical filtration.

In the second part of the paper we study representations of the Drinfel'd quantum double $D(H)$ of algebras $H$ as in the first part. The doubling procedure yields a new class of Hopf algebras beyond a generalization of quantum groups in \cite{{AS1},{AS2}}. For one thing $D(H)$ is not pointed, and for another it does not have decomposition \eqref{triangulatedalgebra}. Nevertheless, $D(H)$ retains enough good features for developing a sort of Lie theory for a Verma-type modules $I(\lambda)$ where $\lambda$ runs over the characters of $G\times\G$. As in part one, but for different reasons, each $I(\lambda)$ has a unique maximal submodule $R(\lambda)$. This enables us to show that the set of simple $D(H)$-modules has the parametrization property. We proceed to describe the radical filtration of $I(\lambda)$ under a mild restriction on the datum for $H$, which is void whenever the certain structure constants $q_i$ of $H$ have odd orders. The `odd order' condition is one way to have all weight spaces of $I(\lambda)$ one-dimensional. When this is the case, the lattice $\Lambda(I(\lambda))$ of $D(H)$-submodules is distributive. This is a very strong property which implies that every submodule of $I(\lambda)$ is a unique sum of some local submodules. Thus the lattice $\Lambda(I(\lambda))$ can be recovered from the partially ordered set $\mathcal J$ of local submodules. We close with a classification of elements of $\mathcal J$.

A more detailed description of material by sections is as follows. In section \ref{prelims} we review the construction of liftings of quantum linear spaces and we develop formulas for skew-derivations assiciated to linked liftings. The main result of the section is the construction of an iterated Ore extension corresponding to the datum of a lifting. This result complements Ore extensions considered in \cite{BDG} and gives an alternate proof of the basis theorem \cite[Prop. 5.2]{AS1}.

In section \ref{H*} we construct the dual basis to the basis of $H$ obtained above. It transpires that the algebra structure of $H^*$ is that of a nilpotent lifting of the quantum linear space with the grouplike and characters switched around. However, $H^*$ is not pointed. Its coradical is computed in \S\ref{coradicalH*}. For related material see \cite{Be}.

In section \ref{triangulated} we first establish a general parametrization theorem for simple $A$-modules where $A$ satisfies \eqref{triangulatedalgebra}-\eqref{nearlynilpotent}. We then turn to nilpotent type liftings and determine the structure of the radical filtration of induced modules $Z(\gamma),\gamma\in\G$, in Theorem \ref{Lseries1}.

We take up the harder case of Drinfel'd double of $H$ in section \ref{D(H)}. Our calculation of multiplication in $D(H)$ is informed by Lemma \ref{xdata} which says that multiplication in the double of a Hopf algebra generated by grouplike and skew-primitive elements is expressed in terms of automorphisms and skew-derivations associated to skew-primitives. In the case at hand we compute explicitly those skew-derivations in a series of lemmas in \S\ref{multiplicationD(H)}. As a first step toward parametrization theorem for simple $D(H)$-modules we find a basic sublagebra in the sense of representation theory of algebras. We then construct a family of induced modules $I(\lambda)$ parametrized by $\lambda\in\GG$ where $\Gamma=G\times\G$. The action of $\Gamma$ splits up $I(\lambda)$ into a direct sum of weight subspaces. These are made explicit in Lemma \ref{weightdecomp} leading up to the parametrization Theorem \ref{simpleD-mod}.

The problem of determining the Loewy filtration of $I(\lambda)$ is finer and it is there that we impose a restriction on datum in Definition \ref{halfcleandata}. The key step of our analyses consists in showing that generators of $D(H)$ act as raising and lowering operators on the weight basis of $I(\lambda)$. From this we derive the Loewy structure of $I(\lambda)$ in Theorem \ref{Lseries2}. The distributive case is handled in Theorem \ref{classicalliftings}.
%%%%%%%%%%%%%%%%%%%%%%%%%%%%%%%%%%%%%%%%%%%%%%%%%%%%%%%%%%%%%%%%%%%%%%%%%%%%%%%%%%%%%%%%%%%%%%%%%%%%%%%%%%%%%%%%%%%%%%%%%%%%%%%%%%%%%%%%%%%%%%%%%%%%%%%%
\section{Preliminaries}\label{prelims}
\subsection{Liftings of $V$}\label{liftings}

We fix some notation. Below $\k$ ia a field of characteristic $0$ containing all roots of $1$ and $\k^{\bullet}=\k\setminus\{0\}$. We denote a finite abelian group by $G$, let $\G:={\rm Hom}(G,\k^{\bullet})$ denote the dual group, and we let $\k G$ stand for the group algebra of $G$ over $\k$. The order of $g\in G$ of a group $G$, is denoted by $|g|$. In particular, for a root of unity $q\in\Bbbk^{\bullet}$, $|q|$ denotes the order of $q$. We set ${\underline n}=\{1,2,\ldots,n\}$ and $[n]:=\{0,1,\ldots,n-1\}$. For a vector space $V$ and a subset $X$ of $V$ we denote the span of $X$ by $(X)$. The unsubscribed `$\otimes$' means `$\otimes_{\Bbbk}$'. For all $n,m\in\mathbb Z$ with $m\ge 0$ ${\binom nm}_q$ denotes the Gaussian $q$-binomial coefficient \cite{KC}, $(n)_q={\binom n1}_q$ and $(n)_q!=(1)_q\cdots (n)_q$.

We review construction of the underlying algebras of this paper following \cite{AS1}-\cite{AS2}. They belong in the class of Hopf algebras parametrized by some elements of $G\times\G$. The starting point of their construction is a left-left Yetter-Drinfel'd finite-dimensional module $V$ over $\Bbbk G$, or a $YD$-module, for short. This means \cite{AS2} that $V$ is a left $\Bbbk G$-module and a left $\Bbbk G$-comodule with the $G$-action preserving $G$-grading. Let's denote by $\omega:V\to\Bbbk G\o V$ the comodule structure map and by `.' the $G$ action. By \cite[Section 1.2]{AS2} $V$ has a basis $\{v_i|i\in\un\}$, where $n=\dim V$, of $G$- and $\G$-eigenvectors, namely there are $a_i\in G$ and $\chi_i\in\G,\,i\in\un$ such that
\begin{align}g.v_i&=\chi_i(g)v_i\label{YD1}\\ 
                                 \omega(v_i)&=a_i\otimes v_i\label{YD2}\end{align}
for all $g\in G$. We set $$\overline{\mathcal D}=(G,(a_i),(\chi_i)|i\in\un)$$ and call this tuple a linear datum associated with $V$. We denote by ${}_G^G{\mathcal YD}$ the category of all $YD$-modules over $\Bbbk G$.

The $YD$-module structure on $V$ extends to a $YD$-structure on $V^{\o m}$ for every integer $m\ge 0$ by using the diagonal action and the codiagonal coaction of $\Bbbk G$ on a tensor product. Explicitly, this means that
\begin{align}g.(v_{i_1}\cdots v_{i_m})&=g.v_{i_1}\cdots g.v_{i_m}.\label{YD1'}\\
\omega(v_{i_1}\cdots v_{i_m})&=a_{i_1}\cdots a_{i_m}\o v_{i_1}\cdots v_{i_m}.\label{YD2'}\end{align}
for all $g$ and $1\le i_1,\ldots,i_m\le n$. Let $F(V)$ be a free associative algebra generated by $V$. As $F(V)=\oplus_{m\ge 0}V^{\o m}$, $F(V)$ becomes a graded $YD$-module and it follows easily from the formulas \eqref{YD1'},\eqref{YD2'} that $F(V)$ is an algebra in ${}_G^G{\mathcal YD}$. Moreover, $F(V)$ is endowed with a special structure of Hopf algebra in ${}_G^G{\mathcal YD}$ \cite [Section 2.1]{AS2}. This is done as follows. First, it is well known \cite{RT,AS2} that the ${}_G^G{\mathcal YD}$ is a braided tensor category with the tensor product just defined and the braiding given by the formula
\begin{equation*}c(u\o v)=u_{(-1)}.v\o u_{(0)},\end{equation*}
where we write $\omega(u)=u_{(-1)}\o u_{(0)}$ for all $u\in F(V)$. Using this braiding we define the multiplication `$\bullet$' in $F(V)\o F(V)$ by
\begin{equation}\label{twisted}(x\o y)\bullet(u\o v)=x(y_{(-1)}.u)\o y_{(0)}v\end{equation}
Second, by a straightforward verification (see also \cite [\S 10.5]{Mo}) the definition \eqref{twisted} turns $F(V)\o F(V)$ into an algebra in  ${}_G^G{\mathcal YD}$ denoted by $F(V)\underline{\o}F(V)$. Since $F(V)$ is a free algebra there is an algebra homomorphism 
$$\delta:F(V)\to F(V)\underline{\o}F(V)$$
defined on the generators by $\delta(v_i)=1\o v_i+v_i\o 1$ for $i\in\un$. Another verification shows that $\delta$ is $G$-linear and $G$-colinear. All in all we see that $F(V)$ is a bialgebra in ${}_G^G{\mathcal YD}$. Further by \cite [11.0.10]{Sw} the coalgebra $F(V)$ has the coradical $\Bbbk$, and then an argument of Takeuchi \cite [5.2.10]{Mo} proves existence of the antipode. Thus $F(V)$ is indeed a Hopf algebra in ${}_G^G{\mathcal YD}$. 

Multiplication law \eqref{twisted} can be elucidated as follows. The formulas \eqref{YD1'}-\eqref{YD2'} allow us to associate with a monomial $\uv=v_{i_1}\cdots v_{i_m}$ the bidegree $(\chi_{\uv},g_{\uv})$ where $\chi_{\uv}=\chi_{i_1}\cdots\chi_{i_m}$ and $g_{\uv}=a_{i_1}\cdots a_{i_m}$. The set $\{\uu\o\uv\}$ forms a basis of $F(V)\o F(V)$ in which the definition \eqref{twisted} takes on the form
\begin{equation}\label{lusztig}(\ux\o \uy)\bullet(\uu\o \uv)=\chi_{\uu}(g_{\uy}){\underline x}\ {\underline u}\o \uy\ \uv\end{equation}
Equation \eqref{lusztig} shows that the definition of $F(V)$ is analogous to Lusztig's definition of algebra ${}'{\mathbf f}$ \cite{Lu1}. Moreover, when $G$ is generated by the $a_i$'s and the mapping $a_i\mapsto\chi_i$ is a homomorphism $G\to\G,\,F(V)$ is exactly the Lusztig's type algebra associated to the bilinear form $(,)\,G\times G\to\Bbbk$ defined on the generators by $(a_i,a_j)=\chi_i(a_j)$, for all $i,j\in\un$.
%For, let $\phi:G\to\G denote the above homomorphism. $\phi$ gives rise to a bilinear form $(,):G\times G\to\Bbbk, (g,h)=\phi(g)(h)$. Since $\phi(a_i)=\chi_i, (a_i,a_j)=\chi_i(a_j)$. Regard $F(V)$ as a $G$-graded vector space by giving a monomial $v=v_{i_1}\cdotsv_{i_m}$ the degree $g_v=a_{i_1}\cdots a_{i_m}$. Then Lusztigs multiplication in $F(V)\o F(V)$ is given by
%$$(x\o y)(u\o v)=(g_u,g_y)xu\o yv.$$ Note that $(g_u,g_y)=\phi(g_u)(g_v)=\chi_u(g_v)$, in agreement with \eqref{lusztig}

We can now define a fundamental object of the theory. Let ${\mathcal F}(V)=F(V)\o \Bbbk G$ be the vector space made into a Hopf algebra by the smash product and smash coproduct constructions. By \cite[Theorem 1]{Ra} ${\mathcal F}(V)$, denoted by $F(V)\#\Bbbk G$, is indeed an ordinary Hopf algebra whose bialgebra structure is descibed by 
\begin{align}(u\#g)(v\#h)&=u(g.v)\#gh\label{smashprod}\\
\Delta(u\#g)&=u^{(1)}\#(u^{(2)})_{(-1)}g\o (u^{(2)})_{(0)}\#g\label{smashcoprod}\end{align}
where we write coproduct of $F(V)$ by $\delta(u)=u^{(1)}\o u^{(2)}$. 

The algebras of interest to us are tied to a special kind of $YD$-module.
\begin{Def}(\cite{AS1}) A $YD$-module $V$ with datum $\overline{\mathcal D}=((a_i),(\chi_i)|i\in\un)$ is called a {\em quantum linear space} if
\begin{equation}\label{qls}\chi_i(a_j)\chi_j(a_i)=1\;\text{for all}\;i\neq j\end{equation}
\end{Def}
From now on we assume that $V$ is a quantum linear space. We let $q_{ij}=\chi_j(a_i)$ for $i\neq j,\,q_i=\chi_i(a_i)$ and $m_i=|q_i|$. 

A {\em datum} (or compatible datum \cite{AS1}) $\mathcal D$ for $V$ is a triple
$$\mathcal D=(\overline{\mathcal D},(\mu_i),(\lambda_{ij}))$$
composed of the linear datum $\overline{\mathcal D}$ of $V$ and two sets of scalars $(\mu_i)_{i\in\un}$ and $(\lambda_{i,j})$ with $i\neq j,\,i,j\in\un$ such that
\begin{align}\mu_i&=0\;\text{if}\;a_i^{m_i}=1\;\text{or}\;\chi_i^{m_i}\neq\epsilon\label{dat1}\\
\lambda_{ij}&=0\;\text{if}\;a_ia_j=1\;\text{or}\;\chi_i\chi_j\neq\epsilon\label{dat2}\\
\lambda_{ji}&=-q_{ji}\lambda_{ij}\label{dat3}\end{align}
We will identify $v_i$ with $v_i\# 1$. For a datum $\mathcal D$ we define the elements $p_i$ and $r_{ij}$ by
\begin{align*}p_i&=v_i^{m_i}-\mu_i(a_i^{m_i}-1)\;\text{for all}\,i\in\un\\
r_{ij}&=v_iv_j-q_{ij}v_jv_i-\lambda_{ij}(a_ia_j-1),\;\text{for all}\,1\le i\neq j\le n\end{align*}
We let $I({\mathcal D})$ be the ideal of ${\mathcal F}(V)$ generated by $p_i, r_{ij}$ for $i,j\in\un$ and we set
$$H({\mathcal D})={\mathcal F}(V)/{I({\mathcal D})}.$$
We remark that formula \eqref{smashcoprod} implies readily that $\Delta(v_i)=v_i\o 1+a_i\o v_i$, thus $S(v_i)=-a_i^{-1}v_i$, where $S$ is the antipode of ${\mathcal F}(V)$. A direct verification yields that $p_i$ is $(a_i^{m_i},1)$-primitive (\cite[p.23]{AS1}) and likewise $r_{ij}$ is $(a_ia_j,1)$-primitive, thanks to \eqref{qls}. In addition, a routine calculation gives $S(p_i)=-a_i^{-m_i}p_i$ \cite[p.24]{AS1} and also $S(r_{ij})=-a_ia_jr_{ij}$. Consequently $I(\mathcal D)$ is a Hopf ideal, hence $H({\mathcal D})$ is a Hopf algebra associated to $\mathcal D$.

As a point of terminology we recall the meaning of lifting of a Hopf algebra \cite{AS2}. A pointed Hopf algebra $H$ is a {\em lifting} of a Hopf algebra $K$ if there is a Hopf algebra isomorphism
\begin{equation}\text{gr}\ H\simeq K,\end{equation}
where $\text{gr}\ H$ is the graded Hopf algebra associated to the coradical filtration of $H$.
Setting the parameters $\mu_i$ and $\lambda_{ij}$ of $\mathcal D$ to zero results in a linear datum $\overline{\mathcal D}$. The Hopf algebra $H(\overline{\mathcal D})$ has a special place in the theory. It is a biproduct of the braided Hopf algebra $R=F(V)/{I(\overline{\mathcal D})}$ and $\Bbbk G$ and by \cite[5.3]{AS1} $\text{gr}\ H(\mathcal D)\simeq H(\overline{\mathcal D})$. As $V$ determines $H(\overline{\mathcal D})$ we call $H(\mathcal D)$ a lifting of $V$. We note that one of the main results of \cite{AS1} is that every lifting of $H(\overline{\mathcal D})$ has the form $H(\mathcal D)$ for a datum $\mathcal D$ on a quantum linear space $V$.
%%%%%%%%%%%%%%%%%%%%%%%%%%%%%%%%%%%%%%%%%%%%%%%%%%%%%%%%%%%%%%%%%%%%%%%%%%%%%%%%%%%%%%%%%%%%%%%%%%%%%%%%%%%%%%%%%%%%%%%%%%%%%%%%%%%%%%%%%%%%%%%%%%%%%%%%
\subsection{Skew-derivations}\label{sderivations}
We begin by recalling the concept of left and right skew-derivation. Let $A$ be an algebra, $a$ an element of $A$ and $\phi$ be an algebra endomorphism of $A$. The assignments
\begin{align*}b\mapsto ab-\phi(b)a,\;\text{for all}\;b\in A\\
              b\mapsto ba-a\phi(b),\;\text{for all}\;b\in A\end{align*}
define two linear mappings denoted by ${}_{\phi}[a,b]$ and $[a,b]_{\phi}$ and called the left and right $\phi$-commutators, respectively. They are a left and right $\phi$-derivations, respectively, in the sence of having the property
\begin{align}{}_{\phi}[a,uv]&={}_{\phi}[a,u]v+\phi(u){}_{\phi}[a,v]\;\text{for all}\;u,v\in A\label{prodrule1}\\
[a,uv]_{\phi}&=[a,u]_{\phi}\phi(v)+u[a,v]_{\phi}\;\text{for all}\;u,v\in A\label{prodrule2}\end{align}
In all applications the endomorphism $\phi$ is the inner automorphism $\imath_g: h\mapsto ghg^{-1},\,h\in H$ induced by an invertible element $g\in H$. We shall use a shorter notation ${}_g[a,b]$ and $[a,b]_g$ for a left/right $\imath_g-\text{commutators}$.

We shall need a commutation formula for powers of generators. Let $A$ be an algebra, $a,b,x,y\in A$ and $\lambda,q\in\Bbbk^{\bullet}$, where $a,b$ are invertible. Suppose
$$gxg^{-1}=qx,\,gyg^{-1}=q^{-1}y\;\text{for}\; g=a,b\;\text{ and}\;{}_b[y,x]=\lambda(ab-1)$$
\begin{Lem}\label{commutatorrel1} For every natural $m\ge 1$ the folowing hold
\begin{align*}{}_{b^m}[y^m,x]&=\lambda(m)_q(q^{m-1}ab-1)y^{m-1}\tag{1}\\
{}_b[y,x^m]&=\lambda(m)_qx^{m-1}(q^{m-1}ab-1)\tag{2}\end{align*}
\end {Lem}
\pf (1) The formula holds for $m=1$ by definition. We induct on $m$ assuming the formula holds for a given $m$. We begin by noting that ${}_{b^m}[y^m,x]=[x,y^m]_{a^{-1}}$. Therefore we can apply \eqref{prodrule2} to carry out the induction step. This gives 
\begin{align*}[x,y^{m+1}]_{a^{-1}}&=y^m[x,y]_{a^{-1}}+[x,y^m]_{a^{-1}}qy\\\intertext{(which by the basis of induction and the induction hypothesis)}&=\lambda[y^m(ab-1)+q(m)_q(q^{m-1}ab-1)y^m].\end{align*}
Since $y^mab=q^{2m}aby^m$ the right hand side equals
\begin{align*}&\lambda(q^{2m}ab-1+q(m)_q(q^{m-1}ab-1))y^m\\&=\lambda(m+1)_q(q^mab-1)y^m\end{align*}
which gives the desired formula.

Part (2) is proven by a similar (and simpler) argument.\qed

We proceed to the general case. The formula below is a generalization of the Kac's formula \cite[(1.3.1)]{DeCK}. 
\begin{Lem}\label{commutatorrel2} For all integers $j$ and $k$ 
$${}_{b^j}[y^j,x^k]=\sum_{i=1}^{\text{min}(j,k)}x^{k-i}f_i^{j,k}y^{j-i}$$
holds, where $\displaystyle f_i^{j,k}=\lambda^i\displaystyle{\binom ji}_q{\binom ki}_q(i)_q!q^{(k-i)(j-i)}\prod_{m=1}^i(q^{j+k-m-i}ab-1)$
\end{Lem}
\pf We can assume $\lambda=1$ by rescaling $x$ via $x'=x/\lambda$ and return back to $x$ by multiplying the formula by $\lambda^k$. The assertion holds for every $j\ge 1$ and $k=1$ by the preceding lemma. We induct on $k$ assuming the lemma holds for every $j\ge 1$ for a given $k$. By \eqref{prodrule2} we carry out the induction step as follows
\begin{align}\label{step1}{}_{b^j}[y^j,x^{k+1}]&={}_{b^j}[y^j,x^k]x+b^jx^kb^{-j}{}_{b^j}[y^j,x]\\
&={}_{b^j}[y^j,x^k]x+q^{kj} x^k{}_{b^j}[y^j,x].\nonumber\end{align}
By the preceding lemma and the induction hypothesis the right hand side of \eqref{step1} equals to
\begin{equation}\label{step2}(\sum_{i=1}^{\text{min}\{j,k\}}x^{k-i}f_i^{j,k}y^{j-i})x+q^{kj}x^k(j)_q(q^{j-1}ab-1)y^{j-1}\end{equation}
We apply Lemma \ref{commutatorrel1} to $y^{j-i}x$ for every $1\le i<j$ 
$$y^{j-i}x=q^{j-i}xy^{j-i}+(j-i)_q(q^{j-i-1}ab-1)y^{j-i-1}$$
and use $(ab)x=q^2x(ab)$ to rewrite $f_i^{j,k}x=x\tilde{f}_i^{j,k}$ where
$$\tilde{f}_i^{j,k}={\binom ji}_q{\binom ki}_q(i)_q!q^{(k-i)(j-i)}\prod_{m=1}^i(q^{j+k+2-m-i}ab-1).$$ 
This allows us to obtain
\begin{align}\label{step3}x^{k-i}f_i^{j,k}y^{j-i}x&=q^{j-i}x^{k+1-i}\tilde{f}_i^{j,k}y^{j-i}\\&+(j-i)_qx^{k+1-(i+1)}f_i^{j,k}(q^{j-i-1}ab-1)y^{j-i-1}\nonumber\end{align}
It follows that 
$${}_{b^j}[y^j,x^{k+1}]=\sum_rx^{k+1-r}f_r^{j,k+1}y^{j-r}$$
with $1\le r\le k+1$, if $k<j$, and $1\le r\le j$, otherwise, thus showing that $1\le r\le\text{min}\{j,k+1\}$. Moreover, by \eqref{step3} $f_r^{j,k+1}$ satisfy the recurrence relation
\begin{align*}\label{recurrencerel}f_1^{j,k+1}&=q^{kj}(j)_q(q^{j-1}ab-1)+q^{j-1}\tilde{f}_1^{j,k}\\
f_r^{j,k+1}&=(j-r+1)_qf_{r-1}^{j,k}(q^{j-r}ab-1)+q^{j-r}\tilde{f}_r^{j,k}\;\text{for all}\;2\le r\le k\\
f_{k+1}^{j,k}&=(j-k)_qf_k^{j,k}(q^{j-k-1}ab-1)\;\text{if}\;k<j\nonumber\end{align*}
We will show that $f_r^{j,k+1},2\le r\le k$ has the desired form leaving verification of the other cases to the reader. To this end we note that $${\binom jr}_q(r)_q!={\binom j{r-1}}_q(r-1)_q!(j-r+1)_q, \;{\binom kr}_q={\binom k{r-1}}_q\frac{(k-r+1)_q}{(r)_q}$$ and
$$\prod_{m=1}^r(q^{j+k+2-r-m}ab-1)=\prod_{m=1}^{r-1}(q^{j+k+1-r-m}ab-1)(q^{j+k+1-r}ab-1).$$ Therefore
$$f_r^{j,k+1}={\binom j{r-1}}_q{\binom k{r-1}}_q(r-1)_q\prod_{m=1}^{r-1}(q^{j+k+1-r-m}ab -1)\phi$$
where $\phi$ can be written in the form $\phi=(j-r+1)_qq^{(k+1-r)(j-r)}\psi$ with 
$$\psi=q^{k+1-r}(q^{j-r}ab-1)+\frac{(k-r+1)_q}{(r)_q}(q^{j+k+1-r}ab-1).$$
It is a straightforward calculation to deduce that\newline $\psi=(q^{k+j+1-2r}+1)\frac{(k+1)_q}{(r)_q}$, which in turn implies that 
$$\phi=(j-r+1)_qq^{(k+1-r)(j-r)}\frac{(k+1)_q}{(r)_q}(q^{k+j+1-2r}ab-1)$$
and this completes the proof.\qed
%%%%%%%%%%%%%%%%%%%%%%%%%%%%%%%%%%%%%%%%%%%%%%%%%%%%%%%%%%%%%%%%%%%%%%%%%%%%%%%%%%%%%%%%%%%%%%%%%%%%%%%%%%%%%%%%%%%%%%%%%%%%%%%%%%%%%%%%%%%%%%%%%%%%%%%%
\subsection{A basis for $H(\D)$}\label{Hbasis}
Our next goal is to give a simple proof of the basis theorem \cite[Prop. 5.2]{AS1}. We will make use of a connection between our algebras and a construction of ring theory known as Ore extension \cite{GW}. In the present setting this connection was studied in \cite{BDG}.
\begin{Thm}\label{basisH} The following set
\begin{equation}\{v_1^{i_1}\cdots v_n^{i_n}g|g\in G,0\le i_j\le m_j-1,j\in\un\}\end{equation}
is a basis of $H({\D})$. This is the {\em standard} basis of $H$
\end{Thm}
\pf For generalities on Ore extensions we refer to \cite{GW} and we adopt its notation. We induct on $n$ starting with $n=1$. We put $R_0=\Bbbk G$ and define an automorphism $\alpha$ of $R_0$ by setting $\alpha(g)=\chi^{-1}(g)g$ and extending it to $R_0$ by linearity. Next we form a left Ore extension $R=\bigoplus_{n\ge 0}R_0x^n$ of $R_0$ with the automorphism $\alpha$ and the $\alpha$-derivation $\delta=0$. Thus $R$ is a free left $R_0$-module with basis $\{x^i|i=0,1\cdots\}$ whose multiplication is generated by the relations
$$x^ix^j=x^{i+j}\;\text{and}\;xr=\alpha(r)x.$$
Let $\D=\{a,\chi,\mu\}$ be a datum on the set $\{1\}$. Let $I=I(\mathcal D)$ be the ideal of $R$ generated by $x^m-\mu(a^m-1)$, where $m=|\chi(a)|$. It is immediate that $R/I$ is a free left and right $R_0$-module with basis $\{gx^i|0\le i<m\}$ and satisfies the algebra defining relations of $H(\D)$. As these relations imply that $H(\D)$ is a span of the set $\{gx^i|0\le i<m,g\in G\}$, the theorem holds for $n=1$.
% The claim about $R/I$ follows from two observations. First $R/I$ is the sum of $x^iR_0,0\le i\le m-1$ because $x^m\in R_0$. Were the sum $\sum_{i=0}^{m-1}x^iR_0$ not direct we would have a relation $r_0+xr_1+\cdotsx^{m-1}r_{m-1}\in I$. That is to say $u:=\sum x^ir_i=f(x^m-\mu(a^m-1))h$ for some $f,h\in R$. Since $\alpha$ is an automorphism, $\text{deg}\,(f(x^m-\mu(a^m-1)h)=\text{deg}(f)+m+\text{h}\ge m$, while $\text{deg}\,(u)<m$, a contradiction.

Let $\check{\un}=\un\setminus \{1\}$ and denote by ${\D}'$ the restriction of $\D$ to $\check{\un}$. Similarly we let $V'=(v_i|i\in\check{\un}),\,{\mathcal F}'={\mathcal F}(V'),\,I'=I({\D}')$ and $H'=H({\D}')\\={\mathcal F}'/{I'}$. We assume that $H'$ is a free span of\\ $\{gv_2^{i_2}\cdots v_n^{i_n}|g\in G,0\le i_j<m_j,\,j\in\check{\un}\}$. We want to show that $H'$ has an automorphism $\alpha$ and a left $\alpha$-derivation $\delta$ such that
\begin{align*}\alpha(gv_2^{i_2}\cdots v_n^{i_n})&=\chi_1^{-1}(g)\chi_2^{i_2}\cdots\chi_n^{i_n}(a_1)gv_2^{i_2}\cdots v_n^{i_n}\;\text{and}\\
\delta(v_j)&=\lambda_{1j}(a_1a_j-1)\,\text{for all}\, j\in\check{\un}\end{align*}
We note that were this true we could form a left Ore extension $R=H'[x;\alpha,\delta]$ and then pass on to ${\overline R}=R/(x^{m_1}-\mu_1(a^{m_1}-1))$. This $\overline R$ satisfies all algebra relations of $H(\D)$ and has the right dimension $|G|m_1\cdots m_n$, which completes the induction step by the argument used for $n=1$.

To prove the statement about $\alpha,\delta$ we introduce $\Bbbk$-algebra $\mathbb F$ freely generated by $\{v_i|i\in\check{\un}\}$ and the set $\{x_g|g\in G\}$. We define the mappings $\alpha,\delta:{\mathbb F}\to{\mathbb F}$ by setting their values on the generators via
\begin{align*}\alpha(x_g)&=\chi_1^{-1}(g)x_g\;\text{and}\;\alpha(v_i)=\chi_i(a_1)v_i\\
\delta(x_g)&=0\;\text{and}\;\delta(v_i)=\lambda_{1i}(x_{a_1a_i}-1)\end{align*}
and extending them to $\mathbb F$ by requiring $\alpha$ and $\delta$ to be an automorphism and left $\alpha$-derivation, respectively. We can form a left Ore extension ${\mathbb F}[x;\alpha,\delta]$. Let $J$ be the ideal of $\mathbb F$ generated by the elements $r_{g,h}=x_gx_h-x_{gh},\,r_{g,i}=x_gv_i-\chi_i(g)v_ix_g$ and the analogs of $p_i$ and $r_{ij}$ in which every apperance of a group element $g$ is replaced by $x_g$. Clearly ${\mathbb F}/J\simeq H'$ as algebras. We claim that $\alpha$ and $\delta$ factor through to $H'$. This boils down to showing that $J$ is invariant under $\alpha$ and $\delta$.

We begin with inclusion $\alpha(J)\subseteq J$. It is trivial to see that $\alpha(r_{g,h})=\chi_1^{-1}(gh)r_{g,h}$ and $\alpha(r_{g,i})=\chi_1^{-1}(g)\chi_i(a_1)r_{g,i}$. Next we have
$$\alpha(p_i)=\chi_i^{m_i}(a_1)v_i^{m_i}-\mu_i(\chi_1^{-1}(a_i^{m_i})x_{a_i^{m_i}}-1).$$
If $\mu_i\neq 0$, then $\chi_i^{m_i}=\epsilon$ by \eqref{dat1} , and, since $\chi_1^{-m_i}(a_i)=\chi_i^{m_i}(a_1)$ by \eqref{qls} we obtain  $\alpha(p_i)=\chi_i^{m_i}(a_1)p_i$. Next, a simple calculation gives $\alpha(r_{ij})=\chi_i\chi_j(a_1)(v_iv_j-q_{ij}v_jv_i)-\lambda_{ij}(\chi_1^{-1}(a_ia_j)x_{a_ia_j}-1)$. If $\lambda_{ij}\neq 0$ then, first, $\chi_i\chi_j=\epsilon$, and, second, $\chi_1^{-1}(a_ia_j)=\chi_i\chi_j(a_1)$ by \eqref{dat2} and \eqref{qls}, respectively. It follows that $\alpha(r_{ij})=\chi_i\chi_j(a_1)r_{ij}$.

Moving on to the inclusion $\delta(J)\subseteq J$ we note that $\delta(r_{g,h})=0=\delta(r_{g,i})$. The first of these equalities is obvious, and the second follows from $\delta(r_{g,i})=(\chi_1^{-1}(g)-\chi_i(g))\lambda_{1i}x_g(x_{a_1a_i}-1)$ together with the fact that $\lambda_{1i}\neq 0$ implies  $\chi_1^{-1}=\chi_i$ by \eqref{dat2}. In case of $p_i$ we have $\delta(p_i)=\delta(v_i^{m_i})={}_{\alpha}[x,v_i^{m_i}]$ in the notation of \S \ref{sderivations}, and the latter is zero $\mod J$ by Lemma \ref{commutatorrel1} (2).

It remains to compute $\delta(r_{ij})$. A direct calculation gives 
\begin{align*}\delta(r_{ij})&=\lambda_{1i}[(x_{a_1a_i}-1)v_j-\chi_j(a_i)\chi_j(a_1)v_j(x_{a_1a_i}-1)]\tag{*}\\&+\lambda_{1j}[\chi_i(a_1)v_i(x_{a_1a_j}-1)-\chi_j(a_i)(x_{a_1a_j}-1)v_i]\end{align*}
Using the relation $x_gv_k\equiv \chi_k(g)v_kx_g \mod J$ we rewrite (*) as follows
$$\delta(r_{ij})\equiv \lambda_{1i}v_j(\chi_j(a_1a_i)-1)+\lambda_{1j}v_i(\chi_j(a_i)-\chi_i(a_1))\mod J$$
The proof is completed by noting that if $\lambda_{1i}\neq 0$, then the equalities $\chi_j(a_1)=\chi_1^{-1}(a_j)=\chi_i(a_j)$ on account of \eqref{qls} and \eqref{dat2}, respectively, give $\chi_j(a_1a_i)=\chi_i(a_j)\chi_j(a_i)=1$ by \eqref{qls} again. Furthermore, if $\lambda_{1j}\neq 0$, then
$\chi_i(a_1)=\chi_1^{-1}(a_i) (\,\text{by}\, \eqref{qls})=\chi_j(a_i)$, as $\chi_1^{-1}=\chi_j$ by  \eqref{dat2}.\qed
% here are details: $\delta(r_{ij})=\delta(v_i)v_j+\alpha(v_i)\delta(v_j)-\chi_j(a_i)(\delta(v_j)v_i+\alpha(v_j)\delta(v_i)$
%$=\lambda_{1i}(a_1a_i-1)v_j+\chi_i(a_1)v_i\lambda_{ij}(a_1a_j-1)-\chi_j(a_i)(\lambda_{1j}(a_1a_j-1)v_i+\chi_j(a_1)v_j\lambda_{1i}(a_1a_i-1))$. This gives (*).  
%%%%%%%%%%%%%%%%%%%%%%%%%%%%%%%%%%%%%%%%%%%%%%%%%%%%%%%%%%%%%%%%%%%%%%%%%%%%%%%%%%%%%%%%%%%%%%%%%%%%%%%%%%%%%%%%%%%%%%%%%%%%%%%%%%%%%%%%%%%%%%%%%%%%%%

\section{Algebra structure of $H^*$}\label{H*}
\subsection{A basis for $H^*$}\label{basisH*}
We begin by fixing some vector notation. For an $n$- tuple $\ui=(i_1,\ldots,i_n)\in{\mathbb {Z}^{\ge 0}}^n$ and any $n$ noncommuting variables $v_1,\ldots, v_n$ we put $v^{\ui}:=v_1^{i_1}\cdots v_n^{i_n}$. We write $\delta_{\ui,\uj}=\delta_{i_1,j_1}\cdots\delta_{i_n,j_n},\, (\ui)!=\prod\limits_{k=1}^{n}(i_k)_{q_k}$, and $\binom \ui\uj=\prod_{k=1}^{n}{\binom {i_k}{j_k}}_{q_k}$. There ${\binom nm}_q$ denotes the Gaussian $q$- binomial coefficient \cite {KC}. We let $u_k$ stand for the $k$th unit vector $(0\cdots 1\cdots 0)(k\text{th}\,1)$. For two vectors $\ui$ and $\uj$ we write $\ui\le \uj$ if $i_k\le j_k$ for all $k\in\un$.

Every $\gamma\in\G$ gives rise to a functional $\widetilde{\gamma}:H\to \k$ defined by 
\begin{equation}\label{dualgroup} \widetilde{\gamma}(v^{\ui}g)=\delta_{\underline 0,\ui}\gamma(g)\end{equation}
The mapping $\gamma\to\widetilde{\gamma}$ is a group embeding $\G\to H^*$, but not a coalgebra map, if the set $\G$ is given the group-like coalgebra structure. Below we identify $\gamma$ with $\widetilde\gamma$ via that embeding.

For every $g\in G$ we associate a minimal idempotent
$$\epsilon_g=\frac{1}{|G|}\sum_{\gamma\in\G}\gamma(g^{-1})\gamma$$
of $\Bbbk\G$. The natural pairing
$$G\times\G\to\k^{\bullet},\;\langle g,\gamma\rangle\mapsto \gamma(g)$$
induces the canonical isomorphism $G\simeq\widehat{\G}$. It follows that the set\newline $\{\epsilon_g|g\in G\}$ forms a basis of $\Bbbk\G$ dual to the standard basis $\{g|g\in G\}$ of $\Bbbk G$.

We will find useful to have a formula for straightening out certain products. For $\um\le\ui$ we define the scalars 

\begin{align*}\phi(\um,\ui)&=\prod_{p=2}^{n}\chi_p^{m_p}(a_1^{i_1-m_1}\cdots a_{p-1}^{i_{p-1}-m_{p-1}})\end{align*}           
\begin{Lem}\label{straightening}In the foregoing notation, for every $g\in G$
\begin{equation}\label{straightprod}v_1^{m_1}a_1^{i_1-m_1}\cdots v_n^{m_n}a_n^{i_n-m_n}g=\phi(\um,\ui)v^{\um}a^{\ui-\um}g\end{equation}
\end{Lem}
\pf The formula follows immediately from relation \eqref{smashprod}.\qed               

We define the functionals $\xi_i,i\in\un$ by the rule
\begin{equation}\xi_k(v^{\ui}g)=\delta_{u_k,\ui}\quad\text{for every}\,g\in G.\end{equation}
\begin{Lem}\label{singlepowers} For every $c< m_k$ 
\begin{itemize}
\item[(i)] $\xi_k^c(v^{\ui}g)=(c)_{q_k}!\delta_{cu_k,\ui}$.
\item[(ii)] $\xi_k^{m_k}=0\quad\text{for all}\;k\in\un$.
\end{itemize}
\end{Lem}
\pf We begin by noting that in view of\newline $(a_k\o v_k)(v_k\o 1)=q_k(v_k\o 1)(a_k\o v_k)$ the quantum binomial formula \cite{KC} gives
\begin{align*}\Delta(v_k^{i_k})&=\sum_{m_k=0}^{i_k}{\binom {i_k}{m_k}}_{q_k}v_k^{m_k}a_k^{i_k-m_k}\o v_k^{i_k-m_k}\end{align*}
It follows from this together with Lemma \ref{straightening} that
\begin{equation}\label{coproduct}\Delta(v^{\ui}g)=\sum_{\um,\ul}\binom{\ui}{\um}\binom{\uj}{\ul}\phi(\um,\ui)v^{\um}a^{\ui-\um}g\o v^{\ui-\um}g.\end{equation}
Now (i) holds for $c=1$ by the definition of $\xi_k$. Assuming it holds for $c$, the induction step is as follows.
\begin{align}\label{powerofxi}&\xi_k^{c+1}(v^{\ui}g)=\langle\xi_k^c\o\xi_k,\Delta(v^{\ui}g)\rangle\\
&=\sum\binom{\ui}{\um}\phi(\um,\ui)\xi_k^c(v^{\um}a^{\ui-\um}g)\xi_k(v^{\ui-\um} g).\nonumber\end{align}
By the induction hypothesis and the basis of induction\newline $\xi_k^c(v^{\um}a^{\ui-\um}g)=(c)_{q_k}!\delta_{cu_k,\um}$,
and $\xi_k(v^{\ui-\um} g)=\delta_{u_k,\ui-\um}\delta_{\uo,\ul}$. 
It follows readily that the nonzero terms in the right side of (\ref{powerofxi}) satisfy $\um=cu_k,\ui-\um=u_k$. Thus $\ui=(c+1)u_k$. Therefore the sum in (\ref{powerofxi}) equals\newline $\binom{(c+1)u_k}{cu_k}\phi(cu_k,(c+1)u_k)$. It remains to note that $\binom{(c+1)u_k}{cu_k}=(c+1)_{q_k}!$ and $\phi(cu_k,(c+1)u_k)=\chi_k^c(a_1^0\cdots a_{k-1}^0)=1$.

(ii) follows from (i) as $(m_k)_{q_k}=0$.\qed

We can give a formula for the dual basis to the standard basis of $H$. For related results see \cite{Be}. 
\begin{Prop}\label{dualbasis}For every $\uc$ there holds
\begin{itemize}
\item[(1)]$\xi^{\uc}(v^{\ui}g)=(\uc)!\delta_{\uc,\ui}.$
\item[(2)] The set
$$\{[(\uc)!]^{-1}\xi^{\uc}\epsilon_g|0\le c_k<m_k\;\text{for all}\;k\in\un\;\text{and}\;g\in G\}$$
is the dual basis to the standard basis of $H$.
\end{itemize}
\end{Prop}
\pf We begin with (1). We induct on $n$, refering to the preceding lemma for the case $n=1$. Assuming the formula holds for all $\uc$ with $c_n=0$, take $\uc$ with $c_n\neq 0$, and set $\uc'=(c_1,\ldots,c_{n-1},0)$. As $\xi^{\uc}=\xi^{\uc'}\xi_n^{c_n}$ formula \eqref{coproduct} gives
$$\xi^{\uc}(v^{\ui}g)=\sum\binom{\ui}{\um}\phi(\um,\ui)\xi^{\uc'}(v^{\um}a^{\ui-\um}g)\xi_n^{c_n}(v^{\ui-\um} g).$$
By the induction hypothesis and Lemma \ref{singlepowers} the sum equals
$$(\uc')!(c_n)_{q_n}\binom{\ui}{\um}\phi(\um,\ui)$$
where $\um=\uc'$ and $\ui-\um=c_nu_n$. Thus $\ui=\uc'+c_nu_n=\uc$ and it is easy to check that $\binom{\ui}{\um}=1 =\phi(\um,\ui)$. This completes the proof of (1). 

We move to (2). Using formula \eqref{coproduct} we compute
$$\xi^{\uc}\epsilon_h(v^{\ui}g)=\sum\binom{\ui}{\um}\phi(\um,\ui)\xi^{\uc}(v^{\um}a^{\ui-\um}g)\epsilon_h(v^{\ui-\um} g).$$
By part (1) and definition \eqref{dualgroup} the value $\upsilon$ of the above sum is\newline 
$\upsilon=(\uc)!\epsilon_h(g)$, provided $\ui=\um=\uc$, and zero, otherwise, as needed.\qed
\begin{Prop}\label{multipH*}\begin{itemize}\item[(1)] For every $\gamma\in\G$ and $1\le k\le n$
$$\gamma\xi_k=\gamma(a_k)\xi_k\gamma$$
\item[(2)]For all $s,t\in\un$ there holds
$$\xi_s\xi_t=\chi_s(a_t)\xi_t\xi_s.$$
\end{itemize}\end{Prop}
\pf To show (1) we compare the values of $\gamma\xi_k(v^{\ui}g)$ and $\xi_k\gamma(v^{\ui}g)$. Using \eqref{coproduct}, the definition of $\xi_k$ and \eqref{dualgroup}, the first scalar equals $\gamma(a_kg)$ while the second is $\gamma(g)$, provided $\ui=u_k$, and zero, otherwise. 

We proceed to a proof of (2). For $s<t$ the preceding proposition gives $\xi_s\xi_t(v^{\ui}g)=\delta_{u_s+u_t,\ui}\delta_{\uo,\uj}$. In the opposite order using \eqref{coproduct} one can compute easily that $\xi_t\xi_s(v^{\ui}g)=\binom{u_s+u_t}{u_t}\phi(u_t,u_s+u_t)$. Noting that 
$\phi(u_t,u_s+u_t)=\chi_t(a_s)$ by the definition of $\phi$ we conclude that $\xi_t\xi_s=\chi_t(a_s)\xi_s\xi_t$. But then (1) holds for $\xi_s\xi_t$ as well, because $\chi_t^{-1}(a_s)=\chi_s(a_t)$ by \eqref{qls}. 

The last two propositions yield a short alternate proof to \cite[Corollary 5.3]{AS1}. 
\begin{Cor}{\rm(\cite{AS1})}\label{coradfilt} Let $H$ be a lifting of a quantum linear space and $H_r$ be the $r$th term of the coradical filtration of $H$. There holds
$$H_r=(v_1^{i_1}\cdots v_n^{i_n}g|\sum i_j\le r\;\text{and}\;g\in G).$$
\end{Cor}
\pf The ideal $J$ of $H^*$ generated by $\xi_i,i\in\un,$ is nilpotent, thanks to Propositions \ref{dualbasis}(2)-\ref{multipH*}. $J$ is the radical of $H^*$ as $H^*/J\simeq\Bbbk \G$ by Proposition \ref{dualbasis}. The assertion follows from \cite[5.2.9]{Mo}.\qed
       
%%%%%%%%%%%%%%%%%%%%%%%%%%%%%%%%%%%%%%%%%%%%%%%%%%%%%%%%%%%%%%%%%%%%%%%%%%%%%%%%%%%%%%%%%%%%%%%%%%%%%%%%%%%%%%%%%%%%%%%%%%%%%%%%%%%%%%%%%%%%%%%%%%%%%%%%

\section{Triangulated Algebras}\label{triangulated}
\subsection{A general theorem}\label{generaltheorem}
In this section we prove a general form of parametrization property for algebras with certain triangular decomposition. Our proof is similar to \cite[Thm. 1]{KR2}. 

We begin with several preliminary remarks. Let $A$ be a Hopf algebra satisfying the conditions set forth in
\eqref{triangulatedalgebra}-\eqref{nearlynilpotent}. Let's call a subalgebra of $A$ normal if it is stable under both adjoint actions of $A_0$. The restriction of counit $\epsilon$ to $A^+$ has kernel $J^+$. Since $\epsilon$ is $A_0$-linear with respect to both adjoint actions, $J^+$ is an a normal subalgebra. The identities $ax=\sum(\text{ad}_{\ell}a_1)(x)a_2$ and $xa=\sum a_1(\text{ad}_ra_2)(x)$ with $a\in A_0,x\in A^+$ show that $A_0A^+=A^+A_0$, hence $A^{\ge 0}:=A_0A^+$ is a subalgebra of $A$. The splitting $A^+=\Bbbk\oplus J^+$ implies that $A^{\ge 0}=A_0\oplus A_0J^+$. Since $J^+$ is a normal subalgebra, $A_0J^+=J^+A_0$ and therefore $A_0J^+$ is a nilpotent ideal of $A^{\ge 0}$. It follows that all simple $A^{\ge 0}$-modules are pullbacks of simple $A_0$-modules along $A^{\ge 0}\to A_0$. For every simple (left) $A_0$-module $V$ we define $A$-module $Z(V)$ by the formula
\begin{equation}\label{inducedmodule} Z(V)=A\otimes_{A^{\ge 0}}V\end{equation}
For an $A$-module $M$ we denote by $M_0$ the socle of its restriction to $A^{\ge 0}$. We need two auxiliary observations.
\begin{Lem}\label{splitting} For every $A_0$-module there holds
$$Z(V)=V\oplus J^-V.$$
\end{Lem}
\pf Condition \eqref{triangulatedalgebra} implies readily the decomposition
$$A=A^{\ge 0}\oplus J^-A^{\ge 0}.$$
Tensoring this direct sum by $V$ over $A^{\ge 0}$ gives the desired formula.\qed
\begin{Lem}\label{radicalZ(V)} For every simple $A_0$-module $V$ the induced module $Z(V)$ has a unique maximal $A$-submodule contained in $J^-V$.
\end{Lem}
\pf By the preceding lemma and since $A_0J^-=J^-A_0$ the subspace $J^-V$ is a maximal $A_0$-submodule of $Z(V)$. Suppose $M$ is a proper $A$-submodule of $Z(V)$ not contained in $J^-V$. Then $M+J^-V=Z(V)$ and since $J^-$ is nilpotent, the argument of the Nakayama's lemma gives $M=Z(V)$, a contradiction. Now set $R$ equal to the sum of all proper $A$-submodules of $Z(V)$.\qed

We denote the maximal submodule of the above lemma by $R(V)$. We define a family of simple $A$-modules by
$$L(V)=Z(V)/R(V).$$
\begin{Thm}\label{generalparametrization} The mapping $V\mapsto L(V)$ sets up a bijection between the isomorphism classes of simple $A_0$-modules and the isomorphism classes of simple $A$-modules.
\end{Thm}
\pf Let $M$ be a simple $A$-module. Select a simple left $A^{\ge 0}$-submodule $V$ of $M$. We observe that an $A^{\ge 0}$-map $\iota_V:V\to Z(V),\\ \iota_V(v)=1\otimes v$ is universal among all $A^{\ge 0}$-maps of $V$ in $A$-modules. Namely, every $f:V\to M$, $M$ is an $A$-module, can be uniquely extended to $f_{*}:V\to M$ satisfying the equality $f_{*}\iota=f$ via $f_{*}(a\otimes v)=af(v)$. It follows that $M\simeq L(V)$ for some simple $A_0$-module $V$. It remains to show that $L(V)\simeq L(U)$ for two simple $A_0$-modules $V$ and $U$ if and only if $V\simeq U$. To this end it suffices to show that $L(V)_0=V$.

Let $\nu:Z(V)\to L(V)$ be the natural epimorphism. Set $\overline {V}=\nu(V)$ and notice that since $\text{Ker}\,\nu=R(V)\subset J^-V$ we have an isomorphism of $A^{\ge 0}$-modules  $\overline{V}\simeq V$ as well as the decomposition $L(V)=\overline{V}\oplus J^-\overline{V}$. Let $\pi$ be the $A_0$-projection of $L(V)$ on $J^-\overline{V}$. Suppose there is a simple $A^{\ge 0}$-submodule $U$ of $L(V)$ distinct from $\overline{V}$. Set $U'=\pi(U)$ and notice that $U'$ is a simple $A_0$-submodule of $J^-\overline{V}$. Evidently $U'\subset U+\overline{V}$, hence $J^+U'=0$. Therefore by simplicity of $L(V)$ we have $L(V)=AU'=U'+J^-U'$. It follows that $L(V)=J^-\overline{V}$ hence $L(V)=J^-L(V)$ forcing $L(V)=0$, a contradiction.\qed

%%%%%%%%%%%%%%%%%%%%%%%%%%%%%%%%%%%%%%%%%%%%%%%%%%%%%%%%%%%%%%%%%%%%%%%%%%%%%%%%%%%%%%%%%%%%%%%%%%%%%%%%%%%%%%%%%%%%%%%%%%%%%%%%%%%%%%%%%%%%%%%%%%%%%%%
\subsection{Representations of $H$}\label{repsofH}
Let $\mathcal D=(G,(a_i),(\chi_i),(\mu_i),(\lambda_{ij}),\,i,j\in\uN)$ be a datum on a quantum linear $N$-dimensional space. Following \cite [Section 4.1]{RS}) we associate to $\mathcal D$ its {\em linking graph} $\Gamma(\mathcal D)$ which is a simple graph with the vertex set $\uN$ and the edge set of all $(i,j)$ such that $\lambda_{ij}\neq 0$. As usual in graph theory the degree of a vertex $i$ is the number $d(i)$ of all $j$ such that $\lambda_{ij}\neq 0$. We say that $\mathcal D$ is {\em simply linked} datum if $d(i)\le 1$ for all $i$. The simplicity condition is not very severe. For by remark \cite [Section 5]{AS1} $d(i)\le 1$ whenever $|q_i|\ge 3$. The vertices of degree zero give rise to generators of $H(\mathcal D)$ lying in the radical of $H(\mathcal D)$. For our purposes we can assume that $\Gamma(\mathcal D)$ does not have such vertices. We call $\mathcal D$ and $H({\mathcal D})$ of {\em nilpotent type} if $\mu_i=0$ for every vertex $i$.
 From now on $\mathcal D$ is a simply linked datum of nilpotent type with every vertex of degree $1$. Clearly, the number of vertices $N$ is even, so we set $n=N/2$. Renumbering verices, if necessary we can assume that the edge set of $\Gamma({\mathcal D})$ is $\{(i,i+n)|i\in\un\}$. It will be convenient to modify notation. We put $b_i=a_{i+n},\,x_i=v_i$ and $y_i=v_{i+n}$ for every $i\in\un$. Rescaling $x_i$ we will assume $\lambda_{i,i+n}=1$. Thus $\mathcal D$ has the form
$$\D=\{G, (a_i),(b_i),(\chi_i), \lambda_{ij},{\underline n}|i,j\in\underline n\}.$$
Now \eqref{dat2} implies $a_ib_i\neq 0$ and $\chi_{i+n}=\chi^{-1}_i$ for all $i\in\un$ which in turn gives the following conditions:
\begin{align*}
a_ib_i&\not= 0\quad \text{for all}\; i\tag{D0}\\
\chi_j(a_i)&=\chi_i(b_j)\quad \text{for all}\; i,j\tag{D1}\\
\chi_i(a_j)\chi_j(a_i)&=1\quad \text{for all}\; i\not= j\tag{D2}\\
\chi_i(b_j)\chi_j(b_i)&=1\quad \text{for all}\; i\not= j\tag{D3}
\end{align*}
The Hopf algebra $H=H(\D)$ attached to $\D$ is explicitly described as follows. $H$ is generated by $G$ and $2n$ symbols $\{x_i,y_i,\,i\in {\underline n}\}$ subject to the relations of $G$ and the following relations:
\begin{itemize}
\item[(R1)] $gx_i=\chi_i(g)x_ig$\quad for all $g\in G$.\\
\item[(R2)] $gy_i=\chi^{-1}(g)y_ig$\quad for all $g\in G$.\\
\item[(R3)] $x_ix_j- q_{ij}x_jx_i=0$\quad for $i\not= j$.\\
\item[(R4)] $y_iy_j-q_{ij}y_jy_i=0$ \quad for $i\not= j$\\
\item[(R5)] $x_iy_j-q^{-1}_{ij}y_jx_i=\delta_{ij}(a_ib_i-1)$\quad for all $i$.\\
\item[(R6)] $x_i^{m_i}=0=y_i^{m_i}$\quad for all $i$.\\
\item[(R7)] $\Delta x_i=a_i\o x_i+ x_i\o 1$\quad for all $i$.\\
\item[(R8)] $\Delta y_i=b_i\o y_i+y_i\o 1$\quad for all $i$. 
\end{itemize}

We proceed now to a classification of simple $H$-modules. Our\newline approach is a generalization of \cite{CK}. Let $Y$ be the subHopfalgebra of $H$ generated by $G$ and $y_i,i\in\un$. For every $\gamma\in\G$ we make $\Bbbk$ a $Y$-module denoted by $\Bbbk_{\gamma}$ by setting
\begin{align*}g.1_{\gamma}&=\gamma(g)\;\text{for all}\;g\in G\\
y_i.1_{\gamma}&=0\;\text{for all}\;i\end{align*}
where $1_{\gamma}$ is identified with $1\in\Bbbk$. We define the $H$-module $Z(\gamma)$ by inducing from $Y$ to $H$, viz.
$$Z(\gamma)=H\otimes_Y\Bbbk_{\gamma}.$$
Since $H$ is a free $Y$ module on a basis $\{x^{\ui}\}$ we see that the set $\{x^{\ui}\otimes 1_{\gamma}\}$ forms a basis for $Z(\gamma)$.

Let $M$ be an $H$-module. We say that $0\neq m\in M$ is a weight element of weight $\gamma\in\G$ if $g.m=\gamma(g)m\;\text{for all }\;g\in G$. A weight element is called primitive if $y_i.m=0\;\text{for all}\;i$. For $\gamma\in\G$ we define $S(\gamma)$ to be the subset of all $j\in\un$ such that 
\begin{equation}\label{special}\gamma(a_jb_j)=q_j^{-e_j}\;\text{for some }\;0\le e_j\le m_j-2.\end{equation}
We denote by $e_j(\gamma)$ the above integer, dropping $\gamma$ whenever it is clear from the context. We say that elements $x,y$ of an algebra skew commute if $xy=qyx$ for some non-zero $q\in\k$.
\begin{Lem}\label{primitivemonomial} A monomial $x^{\ui}\otimes 1_{\gamma}$ is primitive if and only if $i_j=0,e_j+1$ for every $j\in S(\gamma)$, and $i_k=0$ for all $k\notin S(\gamma)$
\end{Lem}
\pf Since $y_j$ skew commutes with every $x_i,i\neq j$, there is $c\in\Bbbk^{\bullet}$ such that 
$$y_j.x^{\ui}\otimes 1_{\gamma}=cx_1^{i_1}\cdots x_{j-1}^{i_{j-1}}y_jx_j^{i_j}\cdots x_n^{i_n}\otimes 1_{\gamma}.$$
By Lemma \ref{commutatorrel1} $y_j.x_j^{i_j}\otimes 1_{\gamma}=-q_j(i_j)_{q_j}x_j^{i_j-1}(q_j^{i_j-1}a_jb_j-1)$, provided $i_j\neq 0$, and $0$, otherwise. Further, for every $k\neq j$ the condition (D1) implies that 
\begin{equation}\label{orthogonality}\chi_k(a_jb_j)=\chi_k(a_j)\chi_k(b_j)=\chi_k(a_j)\chi_j(a_k)=1.\end{equation} 
Therefore $a_jb_j$ commutes with every $x_k,k\neq j$. It follows that $y_j.x^{\ui}\otimes 1_{\gamma}=0$ if and only if $i_j=0,e_j+1$, the last possibility occuring for $j\in S(\gamma)$ only.\qed

By Theorem \ref{generalparametrization} each $Z(\gamma)$ has a unique maximal submodule $R(\gamma)$, possibly zero.  We associate a simple $H$-module 
$$L(\gamma)=Z(\gamma)/{R(\gamma)}$$
to every $\gamma\in\G$. The next result explicitly describes $R(\gamma)$.
\begin{Prop}\label{radical} In the foregoing notation 

{\rm (1)} The family $\{L(\gamma)|\gamma\in\G\}$ is a full set of representatives of simple $H$-modules.

{\rm (2)} $R(\gamma)$ is the sum of all submodules generated by $x_j^{e_j+1}\otimes 1_{\gamma},j\in S(\gamma)$.
\end{Prop}
\pf (1) is a particular case of Theorem \ref{generalparametrization}.

(2) On the one hand each primitive vector $x^{\ui}\otimes 1_{\gamma}$ generates a submodule spanned by all $x^{\uj}\otimes 1_{\gamma}$ with $\uj\ge \ui$, hence a proper one. 

Conversely, suppose $v=\sum_{\ui} c_{\ui}x^{\ui}\otimes 1_{\gamma},c_{\ui}\in\Bbbk^{\bullet}$, generates a proper submodule. If $i_j\ge e_j+1$ for at least one $j\in S(\gamma)$, then $x^{\ui}\otimes 1_{\gamma}$ lies in $H.(x^{e_j+1}\otimes 1_{\gamma})$, hence is contained by $R(\gamma)$ by the opening remark. Suppose $v$ involves a monomial $x^{\ui}\otimes 1_{\gamma}$ with $i_j\le e_j$ for all $j\in S(\gamma)$. By Lemma \ref{commutatorrel2} 
$$y_j^{i_j}x_j^{i_j}\otimes 1_{\gamma}=-q_j(i_j)_{q_j}!\prod_{m=1}^{i_j}(q_j^{i_j-m}\gamma(a_jb_j)-1)\otimes 1_{\gamma}.$$
It follows that 
$$y_1^{i_1}\cdots y_n^{i_n}v=c_{\ui}d\otimes 1_{\gamma}+\sum_{\um\neq \uo}\kappa_{\um}x^{\um}\otimes 1_{\gamma}$$
with $d\neq 0,\kappa_{\um}\in\Bbbk$. But the latter element generates $Z(\gamma)$ since $x_i$ is nilpotent for all $i$. This completes the proof.\qed

We derive from the previous proposition dimension of $L(\gamma)$. 
\begin{Cor}\label{dimensionHsimple} $\dim\,L(\gamma)=\prod_{k\notin S(\gamma)} m_k\prod_{j\in S(\gamma)}(e_j+1)$.\qed
\end{Cor}
We proceed to calculation of the radical and socle series of $Z(\gamma)$. First we recall \cite{LA} that for a finite-dimensional algebra $A$ and a left $A$-module $M$ the radical $R(M)$ of $M$ is the smallest submodule such that $M/{R(M)}$ is semisimple. It is easy to see that $R(M)=JM$ where $J$ is the radical of $A$. Dually the socle $\Sigma(M)$ of $M$ is the largest semisimple submodule of $M$. The radical or Loewy series $\{R^n(M)\}$ of $M$ is defined recursively by $R^0(M)=M$ and $R^m(M)=R(R^{m-1}(M))$ for $m\ge 1$. Similarly, the socle series $\{\Sigma_m(M)\}$ is defined by $\Sigma_0(M)=0$ and $\Sigma_m(M)$ is the preimage in $M$ of $\Sigma(M/{\Sigma_{m-1}(M)})$ for $m\ge 1$. We note that the numbers $\text{min}\,\{m|R^m(M)=0\}$ and $\text{min}\{m|\Sigma_m(M)=M\}$ coincide. The common value is known as the Loewy length of $M$, denoted by $\ell(M)$. Moreover, the two series are related by inclusion $R^m(M)\subseteq\Sigma_{\ell(M)-m}(M)\;\text{for all }\;m$.

For $A=H$ and $M=Z(\gamma)$ we write $\ell(\gamma)$ and $R^m(\gamma),\Sigma_m(\gamma)$ for the Loewy length and the terms of the Loewy and the socle series, respectively. We define {\em the rank} of monomial $x^{\ui}\o1_{\gamma}$ as the number $\text{rk}\,(x^{\ui}\o1_{\gamma})$ of all $j\in S(\gamma)$ such that $i_j\ge e_j(\gamma)+1$.
\begin{Thm}\label {Lseries1} {\rm (1)} For every $m<\ell(\gamma)$ $R^m(\gamma)$ is generated by the primitive vectors of rank $m$.

{\rm(2)} The radical and the socle series coincide.

{\rm (3)} The Loewy layers ${\mathbb L}^m(\gamma):=R^m(\gamma)/{R^{m+1}(\gamma)}$ are given by 
$${\mathbb L}^m\simeq\oplus\{L(\eta)|\eta\,\text{is weight of primitive basis vector of rank}\;m\}.$$

{\rm (4)} $\ell(\gamma)=|S(\gamma)|+1$
\end{Thm}
 \pf (1)-(4). We often drop $\gamma$ when it is clear from the context. We induct on $m$. The assertion holds for $m=0$ as $Z(\gamma)$ is generated by $1\o 1_{\gamma}$. Suppose it is true for $R^m$. Pick a primitive vector $v_{\ui}:=x^{\ui}\o 1_{\gamma}$. Clearly $\eta=\gamma\chi^{\ui}$ is weight of $v_{\ui}$. By the universal property of induced modules there is  an epimorphism $\phi:Z(\eta)\to H.v_{\ui}$ sending $1\o 1_{\eta}\mapsto v_{\ui}$. It follows that $R(H.v_{\ui})=\phi(R(\eta))$ which by the previous proposition equals $\sum H\phi(w)$, where $w$ runs over all primitive monomials of $Z(\eta)$ of rank  $1$. We next compute $S(\eta)$. We have $\eta(a_kb_k)=\gamma\chi_k^{i_k}(a_kb_k)$, because $\chi_l(a_kb_k)=1$ for $l\neq k$, as in the proof of Lemma \ref{primitivemonomial}. Noticing that $i_k=0$ for $k\notin S(\gamma)$ and $i_k=0,e_k(\gamma)+1$, otherwise, we arrive at $\eta(a_kb_k)=\gamma(a_kb_k)$ if $k\notin S(\gamma)$ and $\eta(a_kb_k)=q_k^{-e_k(\gamma)},q_k^{-(m_k-e_k(\gamma)-2)}$, otherwise. It follows that $S(\eta)=S(\gamma)$ with $e_k(\eta)=e_k(\gamma)$ or $e_k(\eta)=m_k-e_k(\gamma)-2$ for all $k\in S(\gamma)$. Furthermore, as $w=x_j^{e_j(\eta)+1}\otimes 1_{\eta}$ for some $j\in S(\gamma)$ we get $\phi(w)=x_j^{e_j(\eta)+1}v_{\ui}$. Therefore if $i_j\neq 0$, then in view of $m_j=e_j(\eta)+e_j(\gamma)+2$ and $x_j^{m_j}=0$, we have $\phi(w)=0$. Otherwise, $e_j(\eta)=e_j(\gamma)$, hence  $\phi(w)=x_j^{e_j(\gamma)+1}v_{\ui}$. It follows that every non-zero $\phi(w)$ has rank $m+1$. Moreover, every primitive monomial of rank $m+1$ has the form $\phi(w)$ for a choice of $v_{\ui}$ and $j$. Noting that $R^{m+1}(\gamma)=\sum R(H.v_{\ui})$ where $v_{\ui}$ runs over all primitive monomials, the induction step is complete. This proves (1) from which (4) is an obvious consequence.

(2) Since $R^{\ell}=\Sigma_0$ we may assume by the reverse induction on $m$ that part (2) holds for all $k>m$. Set $M=Z(\gamma)/R^{m+1}(\gamma)$. By the induction hypotheses $R^{m+1}=\Sigma_{\ell-m-1}$, hence $\Sigma(M)=\Sigma_{\ell-m}/R^{m+1}$. Therefore, were $R^m$ a proper submodule of $\Sigma_{\ell-m}$ there would be a simple $H$-module $L$ of $M$ not contained in $R^m/R^{m+1}$. Let $k$ be the largest integer such that $L\subset R^k/R^{m+1}$. We write $\overline{v_{\ui}}$ for the image of $v_{\ui}$ in $M$ and define $\text{rk}\,(\overline{v_{\ui}})$ as $\text{rk}(v_{\ui})$. We claim that $M$ is the span of the images of monomials. For this is true of $H.v_{\ui}$ for a primitive $v_{\ui}$ because the latter is the span of $x^{\uj}v_{\ui}$ where $x^{\uj}$ runs over the standard basis of $H$. By part (1) same holds for $R^m$ for every $m$, hence for $M$. In fact we have 
$$R^m=(v_{\ui}|\text{rk}\,(v_{\ui})\ge m)$$ 
Let $u$ be a generator of $L$ written as 
\begin{equation*}u=\sum c_{\ui}\overline{v_{\ui}},\,0\neq c_{\ui}\in\Bbbk\tag{*}\end{equation*} 
By the choice of $k$ the sum $u_k$ of terms of (*) of rank $k$ is nonzero. Let's call the number of terms for $u_k$ in the sum (*)  {\em the length} of $u_k$. We pick a generator $u$ with $u_k$ of the smallest length. As $k<m\le \ell-1$ for each $\overline{v_{\ui}}$ of rank $k$ there is $j\in S(\gamma)$ with $i_j<e_j+1$ where $e_j=e_j(\gamma)$. Set $u'=x_j^{e_j+1-i_j}u$ and observe that $u'\neq 0$, because distinct terms $\overline{v_{\ui}}$ remain distinct or zero upon multiplication by $x_j^{e_j+1-i_j}$ and $x_j^{e_j+1-i_j}\overline{v_{\ui}}\neq 0$, since it has rank $k+1\le m$. However, $u_k'$ has length smaller than $u_k$, a contradiction. 

(3) By part (1)
\begin{equation*}{\mathbb L}^m=\sum H.\overline{v_{\ui}}\tag{**}\end{equation*}
where $v_{\ui}$ runs over all primitive basis vectors of rank $m$. For each $\ui$ the set $B=\{x^{\uj}\overline{v_{\ui}}|\text{rk}(x^{\uj}v_{\ui})=m\}$ is a basis of $H.\overline{v_{\ui}}$. The proof follows immediately, once if we show that $\overline{v_{\ui}}$ is the only primitive vector of $H.\overline{v_{\ui}}$ within a scalar multiple.

Suppose $u$ is a primitive vector of $H.\overline{v_{\ui}}$. Write out $u$ in basis $B$
$$u=\sum c_{\uj}x^{\uj}\overline{v_{\ui}},\;0\neq c_{\uj}\in\k.$$
Since $v_{\ui}$ is primitive, the argument of Lemma \ref{primitivemonomial} shows that 
\begin{equation*}y_kx^{\uj}\overline{v_{\ui}}=d_{\uj}x^{\uj-u_k}(q_k^{j_k-1}\eta(a_kb_k)-1)\overline{v_{\ui}},\,d_{\uj}\in\k^{\bullet}\tag{***}\end{equation*}
where $\eta$ is the weight of $v_{\ui}$. Monomials $x^{\uj-u_k}\overline{v_{\ui}}$ are distinct elements of $B$, which implies that $y_ku=0$ if and only if $y_kx^{\uj}\overline{v_{\ui}}=0$. This condition must hold for all $k$ and therefore, by equation (***) it is equivalent  to $x^{\uj}\otimes 1_{\eta}$ is primitive in $Z(\eta)$. From the proof of part (1) we have that for every $k\notin S(\gamma),j_k=0$, and for $k\in S(\gamma)$, either $j_k=0$, or $j_k=m_k-e_k(\gamma)-1, e_k(\gamma)+1$. Assuming $j_k\neq 0$, in the first case $x^{\uj}\overline{v_{\ui}}=0$, and in the second $\text{rk}(x^{\uj}v_{\ui})\ge m+1$, hence $x^{\uj}\overline{v_{\ui}}=0$ again. Thus $\uj=0$, so that $u=c\overline{v_{\ui}}$ for some $c\in\k$. 
It follows that every $H.\overline{v_{\ui}}$ is a simple module and the sum (**) is direct. For otherwise, some primitive $\overline{v_{\ui}}$ would be a linear combination of other primitive monomials of ${\mathbb L}^m$, a contradiction.\qed

%%%%%%%%%%%%%%%%%%%%%%%%%%%%%%%%%%%%%%%%%%%%%%%%%%%%%%%%%%%%%%%%%%%%%%%%%%%%%%%%%%%%%%%%%%%%%%%%%%%%%%%%%%%%%%%%%%%%%%%%%%%%%%%%%%%%%%%%%%%%%%%%%%%%%%%%

\subsection{The coradical of $H^*$}\label{coradicalH*}

We denote by $\lh$ and $\rh$ two standard actions of $H$ and $H^*$ on each other \cite[Chapter 5]{Sw}. For every $\gamma\in\G$ we define subcoalgebra $C(\gamma)$ by $C(\gamma)=H\lh\gamma\rh H$.
\begin{Prop}\label{coradH*}The family $\{C(\gamma)|\gamma\in\G\}$ contains every simple subcoalgebra of $H^*$. Thus
$$\text{corad}\,(H^*)=\sum_{\gamma\in\G}C(\gamma).$$
\end{Prop}
\pf It suffices to show that $H\lh\gamma\simeq L(\gamma)$. To this end we observe $g\lh\gamma=\gamma(g)\gamma$ and $y_k\lh\gamma=0$ for all $k$. The first of these equalities is obvious. For the second we compute 
$$(y_k\lh\gamma)(x^{\ui}y^{\uj}g)=\gamma(x^{\ui}y^{\uj}gy_k)=\chi_k^{-1}(g)\gamma(x^{\ui}y^{\uj}y_kg)=0$$
by the definition of $\gamma$. We see that $\gamma$ is a primitive vector of weight $\gamma$. Therefore $H\lh\gamma$ is the image of $Z(\gamma)$ under $\phi:1\o 1_{\gamma}\mapsto\gamma$. It remains to show that $\phi(R(\gamma))=0$. By Proposition \ref{radical} this is equivalent to the equality $x_k^{e_k+1}\lh\gamma=0$ for every $k\in S(\gamma)$. 

Pick an integer $m$. By definition of the left action

\noindent$\upsilon=(x_k^m\lh\gamma)(x^{\ui}y^{\uj}g)=\gamma(x^{\ui}y^{\uj}gx_k^m)$. Since every $g\in G$ and every $y_j,j\neq k$, skew commute with $x_k$ we can reduce $\upsilon$ to the form 
$$\upsilon=c\chi_k^m(g)\gamma(x^{\ui}y^{\uj'}(y_k^{j_k}x_k^m)y^{\uj''}),\;c\in\k^{\bullet},$$
where $\uj'=(j-1,\ldots,j_{k-1}),\uj''=(j_{k+1},\ldots,j_n)$. Using Lemma \ref{commutatorrel2} we see readily that $\upsilon=0$ unless $\ui=\uo,\,\uj=mu_k$. When these conditions hold $\upsilon=\chi_k^m(g)\gamma(f_m^{m,m})\gamma(g)$, where $f_m^{m,m}=(m)_{q_k}!\prod_{p=1}^m(q_k^{m-p}a_kb_k-1)$ again by Lemma \ref{commutatorrel2}. As $\prod_{p=1}^m(q^{m-p}\gamma(a_kb_k)-1)=0$ for every $k\in S(\gamma)$ and $m\ge e_k+1$, we conclude that every $C(\gamma)$ is simple coalgebra.  

On the other hand every simple $H$-module is isomorphic to $L(\gamma)$ by Proposition \ref{radical}, which completes the proof.\qed

The functions 
\begin{equation}\label{functionsofgamma}c_k^m:\G\to\k,\;c_k^m(\gamma)=\prod_{p=1}^m(q_k^{m-p}\gamma(a_kb_k)-1)\end{equation}
will play a r\^ole below.
%%%%%%%%%%%%%%%%%%%%%%%%%%%%%%%%%%%%%%%%%%%%%%%%%%%%%%%%%%%%%%%%%%%%%%%%%%%%%%%%%%%%%%%%%%%%%%%%%%%%%%%%%%%%%%%%%%%%%%%%%%%%%%%%%%%%%%%%%%%%%%%%%%%%%%%%

\section{The Drinfel'd double}\label{D(H)}
\subsection{Multiplication in $D(H)$}\label{multiplicationD(H)}
The original definition of the Drinfel'd double $D(H)$ \cite{Dr} of a Hopf algebra is rather technical. For an intrinsic definition of $D(H)$ via the double crossproduct construction see \cite{Mj1}-\cite{Mj2}. We will follow, though, a more transparent description of $D(H)$ due to Doi-Takeuchi \cite{DT}.

We recall that $D(H)$ is $H^*\o H$ as a vector space and $H^{*\text{cop}}\o H$ as a coalgebra with the tensor product coalgebra structure. Note that if $S$ is the antipode of $H$, then $(S^{-1})^*$ is the antipode of $H^{*\text{cop}}$. There is a natural bilinear form  
$$\tau:H^{*\text{cop}}\o H\to\k,\;\tau(\alpha,h)=\alpha(h),\;\text{for}\;\alpha\in H^*,h\in H.$$ 
$\tau$ is an invertible bilinear form in the convolution algebra\newline $\text{Hom}_{\k}(H^{*\text{cop}}\o H,\k)$ with the inverse $\tau^{-1}(\alpha,h)=\tau((S^{-1})^*(\alpha),h)$. Using $\tau$ the algebra structure on $D(H)$ is given with product
\begin{equation}\label{multD(H)}(\alpha\o h)(\beta\o k)=\alpha\tau(\beta_3,h_1)\beta_2h_2\tau^{-1}(\beta_1,h_3)k\end{equation}
where $\Delta_{H^*}^{(2)}(\beta)=\beta_1\o\beta_2\o\beta_3$ and $\Delta_H^{(2)}(h)=h_1\o h_2\o h_3$. In what follows we will drop the ``$\o$''-sign and write $\alpha h$. The essential part of definition \eqref{multD(H)} is
\begin{align}\label{shortmultD(H)}h\beta&=\tau(\beta_3,h_1)\beta_2h_2\tau^{-1}(\beta_1,h_3)\\
&=\beta_3(h_1)\beta_2h_2\beta_1(S^{-1}(h_3))\nonumber\end{align}
Inverting \eqref{shortmultD(H)} gives the equivalent identity
\begin{align}\label{shortmultD(H)'}\beta h&=\tau^{-1}(\beta_3,h_1)h_2\beta_2\tau(\beta_1,h_3)\\
&=\beta_3(S^{-1}(h_1)h_2\beta_2\beta_1(h_3)\nonumber\end{align}
It is convenient  to rewrite identities \eqref{shortmultD(H)} and \eqref{shortmultD(H)'} in terms of actions $\lh$ and $\rh$. An immediate verification gives
\begin{align}\label{standardmult}h\beta&=(h_1\lh\beta\rh S^{-1}(h_3))h_2\\
&=\beta_2((S^{-1})^*(\beta_1)\lh h\rh\beta_3)\nonumber\quad\text{and}\end{align}
\begin{align}\label{standardmult'}\beta h&=h_2(S^{-1}(h_1)\lh\beta\rh h_3)\\
&=(\beta_1\lh h\rh (S^{-1})^*(\beta_3))\beta_2\nonumber\end {align}
We note that formulas \eqref{standardmult} and \eqref{standardmult'} were obtained in \cite{Ra1} and \cite{Ro}, respectively.

One consequence of \eqref{standardmult} is the formula $g\alpha g^{-1}=g\lh \alpha\rh g^{-1}$. It shows that $H^*$ is invariant under the action by $G$ by conjugation. Therefore we have a Hopf subalgebra $\widetilde {H^*}:=H^*\#\k G$ in $D(H)$. Next, suppose $x$ is a $(a,1)$-primitive element of $H$ satisfying $gx=\chi_x(g)xg$ for $\chi_x\in\G$. We associate with $x$ two mappings $\phi_x,\delta_x:\widetilde{H^*}\to\widetilde{H^*}$ defined as follows
\begin{align*}\phi_x(\alpha g)&=(a^{-1}\lh\alpha)\chi_x(g)g\\
\delta_x(\alpha)&=(\alpha\rh xa^{-1})a-xa^{-1}\lh\alpha\quad\text{and}\\
\delta_x(\alpha g)&=\delta_x(\alpha)\chi_x(g)\end{align*}
for all $\alpha\in H^*$ and $g\in G$.
\begin{Lem}\label{xdata} {\rm (1)}  $\phi_x$ is algebra automorphism and $\delta_x$ is a right $\phi_x$-derivation of $\widetilde{H^*}$.

${\rm (2)}\quad [x^s,\alpha]_{\phi_x^s}=[x^{s-1},\alpha]_{\phi_x^{s-1}}x+x^{s-1}[x,\phi^{s-1}_x(\alpha)]_{\phi}$
for every $\alpha\in H^*$ and $s\ge 1$
\end{Lem}
\pf (1) The claim about $\phi_x$ is obvious. For the rest it suffices to show that
$$\alpha x=x\phi_x(\alpha)+\delta_x(\alpha)\quad\text{for all}\;\alpha\in H^*.$$
Note the equalities $\Delta^{(2)}(x)=a\o a\o x+a\o x\o 1+x\o 1\o 1$ and $S^{-1}(x)=-xa^{-1}$. Therefore we have from \eqref{standardmult'} 
$$\alpha x=a(a^{-1}\lh\alpha\rh x)+x(a^{-1}\lh\alpha)-xa^{-1}\lh\alpha.$$
As  $a(a^{-1}\lh\alpha\rh x)=(\alpha\rh xa^{-1})a$ by \eqref{standardmult} and $gx=\chi_x(g)xg$, the proof of (1) is complete.

(2) We recall that $[x^s,\alpha]_{\phi_x^s}$ stands for the right derivation $\alpha\mapsto \alpha x^s-x^s\phi_x^s(\alpha)$ as defined in Section \ref{sderivations}. The proof of the identity is by a direct verification.\qed
%%%%%%%%%%%%%%%%%%%%%%%%%%%%%%%%%%%%%%%%%%%%%%%%%%%%%%%%%%%%%%%%%%%%%%%%%%%%%%%%%%%%%%%%%%%%%%%%%%%%%%%%%%%%%%%%%%%%%%%%%%%%%%%%%%%%%%%%%%%%%%%%%%%%%%%%
%%%%%%%%%%%%%%%%%%%%%%%%%%%%%%%%%%%%%%%%%%%%%%%%%%%%%%%%%%%%%%%%%%%%%%%%%%%%%%%%%%%%%%%%%%%%%%%%%%%%%%%%%%%%%%%%%%%%%%%%%%%%%%%%%%%%%%%%%%%%%%%%%%%%%%%%

For the sequel we must modify the standard basis of $H$ and the generators of $H^*$. Since every $g\in\G$ skew commutes with all $x_i$ and $y_i$ the set $\{x^{\ui}gy^{\uj}|0\le i_k,j_k<m_k\;\text{and}\;g\in G\}$ is another basis of $H$. We define the functionals $\gamma,\xi_k,\eta_k$ for all $k\in\un$ by setting
\begin{align}{\gamma}(x^{\ui}gy^{\uj})&=\delta_{\underline 0,\ui}\delta_{\underline 0,\uj}\gamma(g)\\
\xi_k(x^{\ui}gy^{\uj})&=\delta_{u_k,\ui}\delta_{\uo,\uj}\quad\text{for every}\,g\in G.\\
              \eta_k(x^{\ui}gy^{\uj})&=\delta_{\uo,\ui}\delta_{u_k,\uj}\quad\text{for every}\,g\in G.\end{align}
By an argument almost identical to one for Lemma \ref{singlepowers} one can show
\begin{Lem}\label{singlepowers'} The formulas \begin{align*}\xi_k^c\gamma(x^{\ui}gy^{\uj})&=(c)_{q_k}!\delta_{\ui,cu_k}\delta_{\uj,\uo}\gamma(g)\;\text{and}\\
\eta_k^c\gamma(x^{\ui}gy^{\uj})&=(c)_{q_k}!\delta_{\ui,\uo}\delta_{\uj,cu_k}\gamma(g).\end{align*}
hold for all $\gamma\in\G,k\in\un$ and $0\le c\le m_k$.
\end{Lem}\qed

We also record an analog of Proposition \ref{multipH*}.
\begin{Prop}\label{multipH*'}\begin{itemize}\item[(1)] For every $\gamma\in\G$ and $1\le k\le n$
$$\gamma\xi_k=\gamma(a_k)\xi_k\gamma\;\text{and}\;\gamma\eta_k=\gamma(b_k)\eta_k\gamma.$$
\item[(2)]For all $s,t\in\un$ the equalities
\begin{align*}\xi_s\xi_t&=\chi_s(a_t)\xi_t\xi_s,\\
\eta_s\eta_t&=\chi_s(b_t)\eta_t\eta_s,\\
\xi_s\eta_t&=\eta_t\xi_s\end{align*}
hold 
\end{itemize}
\end{Prop}\qed

We move on to an explicit description of multiplication in $D(H)$. We start off with the conjugation action of $G$.
\begin{Lem}\label{gconjug} For all $g\in G,\gamma\in\G,1\le k\le n$ the identities 
\begin{align}g\gamma&=\gamma g\\
g\xi_k&=\chi^{-1}_k(g)\xi_kg\label{gxi}\\
g\eta_k&=\chi_k(g)\eta_kg\label{geta}\end{align}
hold.
\end{Lem}
\pf By \eqref{shortmultD(H)} $g\alpha g^{-1}=g\lh\alpha\rh g^{-1}$. From the definition of $\gamma,\xi_k,\eta_k$ the equations
\begin{align}g\lh\gamma&=\gamma(g)\gamma\;\text{and}\;\gamma\rh g=\gamma(g)\gamma\label{ghitgamma}\\
g\lh\xi_k&=\xi_k\;\text{and}\;\xi_k\rh g=\chi_k(g)\xi_k\label{ghitxi}\\
g\lh\eta_k&=\chi_k(g)\eta_k\;\text{and}\;\eta_k\rh g\label{ghiteta}=\eta_k\end{align}
follow which complete the proof.\qed

We need a technical lemma. Below we use the convention that for any set of variables $v_j$, $v^{\ui}=0$ if $i_k<0$ for at least one $k$.
\begin{Lem} There are scalars $c,c',d,d'\in\k^{\bullet}$ depending on $k,\ui,\uj$ and $g,h\in G$ such that
\begin{align}(x^{\ui}gy^{\uj})(x_kh)&=cx^{\ui+u_k}ghy^{\uj}+c'x^{\ui}gh(q_k^{j_k-1}a_kb_k-1)y^{\uj-u_k}\label{rightprod}\\
(y_kh)(x^{\ui}gy^{\uj})&=dx^{\ui}ghy^{\uj+u_k}+d'x^{\ui-u_k}gh(q_k^{i_k-1}a_kb_k-1)y^{\uj}\label{leftprod}\end{align}
\end{Lem}
\pf  Since elements of $G$ and $y_l,l\neq k,$ skew commute with $x_k$, $(x^{\ui}gy^{\uj})(x_kh)=ax^{\ui}ghy^{\uj'}y_k^{j_k}x_ky^{\uj''}$ with $\uj'=(j_1,\ldots,j_{k-1})$ and $\uj''=(j_{k+1},\ldots,j_n)$. Next we use Lemma \ref{commutatorrel1}(1) according to which
$$y_k^{j_k}x_k=q_k^{j_k}x_ky_k^{j_k}-q_k(j_k)_{q_k}(q_k^{j_k-1}a_kb_k-1)y_k^{j_k-1}.$$
This formula and the fact that $x_k$ skew commutes with $x_l,l\neq k,$ completes the proof of \eqref{rightprod}.

The proof of \eqref{leftprod} is almost identical. One must use Lemma \ref{commutatorrel1}(2) together with the observation that $x_l,l\neq k,$ commutes with $a_kb_k$.\qed

The next three lemmas completely determine the algebra structure of $D(H)$.
\begin{Lem}\label{gammacommutators} For every $\gamma\in\G$ and $1\le k\le n$
\begin{align}\gamma x_k&=\gamma(a_k^{-1})x_k\gamma+\gamma(a_k^{-1})q_k(\gamma(a_kb_k)-1)\eta_k\gamma\label{gammax}\\
\gamma y_k&=\gamma(b_k^{-1})y_k\gamma-\gamma(b_k^{-1})(\gamma(a_kb_k)-1)\xi_k\gamma b_k\label{gammay}\end{align}
\end{Lem}
\pf By Lemma \ref{xdata} we need only compute $\delta_{x_k}(\gamma)$. By definition this involves finding $\gamma\rh x_ka_k^{-1}$ and $x_ka_k^{-1}\lh\gamma$. First we show that
\begin{equation}\label{xrhitgamma}\gamma\rh x_ka_k^{-1}=0\end{equation}
For by definition $(\gamma\rh x_ka_k^{-1})(x^{\ui} gy^{\uj})=\gamma(x_ka_k^{-1}x^{\ui}gy^{\uj})$ and the latter is zero because $x_ka_k^{-1}x^{\ui}gy^{\uj}=cx^{\ui+u_k}ga_k^{-1}y^{\uj},c\in\k^{\bullet},$ with $\ui+u_k\neq \uo$ for all $\ui,\uj,g$.

Next we compute 
$$\upsilon:=(x_ka_k^{-1}\lh\gamma)(x^{\ui}gy^{\uj})=\gamma(x^{\ui}gy^{\uj}x_ka_k^{-1}).$$
Using \eqref{rightprod} we express $\upsilon=\upsilon_1+\upsilon_2$ where $\upsilon_1=c\gamma(x^{\ui+u_k}ga_k^{-1}y^{\uj})$ and\linebreak $\upsilon_2=c'\gamma(x^{\ui}ga_k^{-1}(q_k^{j_k-1}a_kb_k-1)y^{\uj-u_k})$. As in the proof of \eqref{xrhitgamma} $\upsilon_1=0$ for all basis elements. The definition of $\gamma$ makes it clear that $\upsilon_2=0$, unless $\ui=\uo$ and $\uj=u_k$. In the latter case $c'=-q_k$, hence $\upsilon_2=-q_k\gamma(ga_k^{-1}(a_kb_k-1))$. That is to say 
\begin{equation}\label{lhitgammax}x_ka_k^{-1}\lh\gamma=-q_k\gamma(a_k^{-1})\gamma(a_kb_k-1)\eta_k\gamma\end{equation} 
by Lemma \ref{singlepowers'}, and this completes the proof of \eqref{gammax}.

The proof of \eqref{gammay} is almost identical. The main steps are the equalities
\begin{equation}y_kb_k^{-1}\lh\gamma=0\;\text{and}\;\gamma\rh y_kb_k^{-1}=-\gamma(b_k^{-1})\gamma(a_kb_k-1)\xi_k\gamma\label{yhitgamma}.\qed\end{equation}

\begin{Lem}\label{xicommutators} For all $k,l\in\un$
\begin{align}\xi_lx_k&=x_k\xi_l+\delta_{k,l}(a_l-\chi_k)\label{xicommx}\\
\xi_ly_k&=y_k\xi_l-\delta_{k,l}q_k^{-1}(q_k-1)\xi_k^2b_k\label{xicommy}\end{align}
\end{Lem}
\pf We use again Lemma \ref{xdata}. Since $a_k^{-1}\lh\xi_l=\xi_l$ by \eqref{gxi}, it remains to find $\delta_{x_k}(\xi_l)$. First, we claim that
\begin{equation}\label{xrhitxi}\xi_l\rh x_ka_k^{-1}=\delta_{k,l}\epsilon\end{equation}
This is the matter of showing that $\upsilon=\xi_l(x_ka_k^{-1}x^{\ui}gy^{uj})=\delta_{k,l}\delta_{\ui,\uo}\delta_{\uj,\uo}$. We observe that $\upsilon=c\chi^{\ui-\uj}(a_k^{-1})\xi_l(x^{\ui+u_k}hy^{\uj})$ for some $h\in G$ and $0\neq c\in\Bbbk$ with $c=1$ when $\ui=\uo=\uj$, which yields the claim by the definition of $\xi_l$.

Next we show the identity
\begin{equation}\label{xlhitxi}x_ka_k^{-1}\lh\xi_l=\delta_{k,l}\chi_k\end{equation}
Now we must compute $\upsilon'=\xi_l(x^{\ui}gy^{\uj}x_ka_k^{-1})$. Applying the straightening out formula \eqref{rightprod} we write $\upsilon'=\upsilon_1+\upsilon_2$ where $\upsilon_1=c\xi_l(x^{\ui+u_k}ga_k^{-1}y^{\uj})$ and $\upsilon_2=c'\xi_l(x^{\ui}ga_k^{-1}(q_k^{j_k-1}a_kb_k-1)y^{\uj-u_k})$. We note that $\upsilon_2$ is zero for all $\ui,g$ and $\uj$. For, if $\uj-u_k\neq \uo$, then surely $\upsilon_2=0$. Else, $\upsilon_2=c'\xi_l(x^{\ui}(s-t))$ for some $s,t\in G$, which is again zero.

As for $\upsilon_1$, if $l\neq k$, then $\upsilon_1=0$, because $\ui+u_k\neq u_l$ for all $\ui$. Suppose $l=k$. Then $\upsilon_1\neq 0$ if and only if $\ui=\uo=\uj$. If so, $c=\chi_k(g)$, hence $\upsilon_1=\chi_k(g)\delta_{\ui,\uo}\delta_{\uj,\uo}$, which gives \eqref{xlhitxi}.

For the proof of the second part we need two observations. First off, the equality $y_kb_k^{-1}\lh\xi_l=0$ for all $k,l$ is self-evident. In the second place we claim the identity
\begin{equation}\label{yrhitxi}\xi_l\rh y_kb_k^{-1}=-q_k^{-1}(q_k-1)\xi_k^2\end{equation}
Again the proof boils down to finding $\upsilon=\xi_l(y_kb_k^{-1}x^{\ui}gy^{\uj})$ which by  \eqref{leftprod} splits up as $\upsilon=\upsilon_1+\upsilon_2$ with $\upsilon_1=d\xi_l(x^{\ui}gb_k^{-1}y^{\uj+u_k})$ and $\upsilon_2=d'\xi_l(x^{\ui-u_k}gb_k^{-1}(q_k^{i_k-1}a_kb_k-1)y^{\uj})$. Now $\upsilon_1$ is always zero, because $\uj+u_k\neq u_l$ for all $\uj$. If $l\neq k$, then $\ui-u_k\neq u_l$ for all $\ui$ forces $v_2=0$ for all choices of $\ui,g,\uj$. Take the case $k=l$. Now $\upsilon_2=0$ for all $\ui,\uj$ such that $\ui-u_k\neq u_k$ or $\uj\neq \uo$. In the remaining case, i.e. $\ui=2u_k$ and $\uj=\uo$, we have from Lemma \ref{commutatorrel1}(2) the identity 
$$y_kb_k^{-1}x_k^2=q_k^{-2}y_kx_k^2b_k^{-1}=(x_k^2y_k-q_k^{-1}(2)_{q_k}x_k(q_ka_kb_k-1))b_k^{-1}$$
which gives $d'=-q_k^{-1}(2)_{q_k}$. Noting that $\xi_k(x_k(q_ks-t))=q_k-1$ for all $s,t\in G$ we arrive at $\upsilon_2=-q_k^{-1}(2)_{q_k}(q_k-1)\delta_{\ui,2u_k}\delta_{\uj,\uo}$. In view of Lemma \ref{singlepowers'} we obtain the desired formula.\qed

\begin{Lem}\label{etacommutators} For all $k,l\in\un$
\begin{align}\eta_lx_k&=q_{kl}^{-1}x_k\eta_l+\delta_{k,l}(q_k-1)\eta_k^2\label{etacommx}\\
\eta_ly_k&=q_{lk}^{-1}y_k\eta_l+\delta_{k,l}q_k^{-1}(\chi_kb_k-\epsilon)\label{etacommy}\end{align}
\end{Lem}
\pf The proof follows from the following equations
\begin{equation}\label{xrlhiteta}\eta_l\rh x_ka_k^{-1}=0\;\text{for all}\;l,k,\;\text{and}\;x_ka_k\lh\eta_l=\delta_{k,l}(q_k-1)\eta_k^2\end{equation}
\begin{equation}\label{yrlhiteta}y_kb_k^{-1}\lh\eta_l=\delta_{k,l}q_k^{-1}\epsilon\;\text{and}\;\eta_l\rh y_kb_k^{-1}=\delta_{k,l}q_k^{-1}\chi_k\end{equation}
A verification of these equations follows the proof of the preceding lemma and is left for the reader.\qed

In keeping with our convention we put $\um=(m_1,\ldots,m_n)$. We let $\mathbb Z(\um)$ denote $[m_1]\times\cdots \times[m_n]$. As a consequence of the last three lemmas we have the fact that the set of all products of $x^{\ui},y^{\uj},\xi^{\uk},\eta^{\ul}$ with $\ui,\uj,\uk,\ul\in\mathbb Z(\um)$ in any prescribed order forms a basis for $D(H)$.

\subsection{Parametrization of Simple $D(H)$-Modules}\label{parametrization}
We denote by $\Gamma$ the group $G\times\G$. We recall that $J(A)$ denotes the radical of algebra $A$. We consider the subalgebra $A=A(\Gamma,y,\xi)$ of $D(H)$ generated by $\Gamma,y_i,\xi_i$ for all $i\in\un$. Relations \eqref{gammay} and \eqref{xicommy} imply that $A$ is a free span of the set $\{y^{\ui}\xi^{\uj}g\gamma|0\le i_k,j_k\le m_k-1, g\gamma\in\Gamma\}$. We set $|\ui|=\sum i_k$ for an $n$-tuple $\ui$.
\begin{Lem}\label{algebraA} $A$ is a right coideal of $D(H)$ and a basic algebra in the sence that
$$A=\k\Gamma\oplus J(A).$$
\end{Lem}
\pf We let $B$ denote the subalgebra of $A$ generated by $\Gamma$ and all $\xi_k$. Clearly $B=\k\Gamma\oplus J(B)$, where $J(B)$ is the span of all $\xi^{\uj}g\gamma$ with $|\uj|>0$. Since the elements of $\Gamma$ skew-commute with every $\xi_k$, $J(B)$ is nilpotent. Let $N$ be the largesr integer with $J(B)^N\neq 0$. Set $I$ to be the span of all $y^{\ui}\xi^{\uj}g\gamma$ with $|\ui|+|\uj|>0$. Evidently $I$ is a complement of $\k\Gamma$ in $A$, and $I=\sum_{|\ui|>0}y^{\ui}B$. We lift the radical filtration 
$$B\supset J(B)\supset\cdots\supset J^N(B)\supset 0$$
to a filtration 
$$I=I_0\supset I_1\supset\cdots\supset I_N\supset 0$$
where $I_k$ is defined by $I_k=\sum_{|\ui|>0}y^{\ui}J^k(B)$. Thanks to the relations \eqref{gammay} and \eqref{xicommy} the elements of $B$ skew commute with all $y^{\ui}z,\,z\in J^k(B)$ modulo $I_{k+1}$. Therefore the $I_k$ form an ideal filtration with nilpotent quotients $I_k/I_{k+1}$ because the $y_k$ generate a nilpotent subalgebra. Thus $I=J(A)$.

We take up the first claim. The comultiplication of $y^{\ui}g$ (cf. \eqref{coproduct}) makes it clear that it suffices to show that the subalgebra $B'$ of $H^*$ generated by $\G$ and the $\xi_k$ is a right coideal. Since $\Delta_D=\Delta_{H^{*\text{cop}}}$ on $H^*$ this is equivalent to $B'$ is a left coideal of $H^*$, or a right $H$-submodule with respect to the `$\rh$'-action. As $H^*$ is an $H$-module algebra we need only to establish inclusion $z\rh h\in B'$ for $z$ and $h$ running over the generators of $B'$ and $H$, respectively. Now using \eqref{ghitgamma}, \eqref{ghitxi}, and \eqref{xrhitgamma}, \eqref{yhitgamma} we have
\begin{align*}\gamma\rh g &=\gamma(g)\gamma,\,\xi_k\rh g=\chi_k(g)\xi_k\;\text{for all}\;g\in G,\quad\text{and}\\
\gamma\rh x_k &=0,\;\text{and}\;\gamma\rh y_k=-(\gamma(a_kb_k)-1)\xi_k\gamma.\end{align*}
A reference to \eqref{xrhitxi} and \eqref{yrhitxi} completes the proof.\qed

The above lemma makes it obvious that every simple $A$-module is $1$-dimensional. Pick $\lambda\in\GG$ and let $\k_{\lambda}$ denote the $A$-module $\k$ with the $A$-action defined by 
$$J(A).1_{\k}=0\;\text{and}\;g\gamma.1_{\k}=\lambda(g\gamma).$$
We write $1_{\lambda}$ for the element $1_{\k}$ of $\k_{\lambda}$. We define a family of $D$-modules $I(\lambda),\lambda\in\GG$ by setting
$$I(\lambda)=D\otimes_A\,\k_{\lambda}.$$
where we write $D=D(H)$. As we mentioned earlier the set\newline $\{x^{\ui}\eta^{\uj}y^{\uk}\xi^{\ul}|\ui,\uj,\uk,\ul\in\mathbb Z(\um)\}$ is a basis for $D(H)$. Therefore, directly from the definition we obtain that $I(\lambda)$ is a free span of the set\newline $\{x^{\ui}\eta^{\uj}\otimes 1_{\lambda}|\ui,\uj\in\mathbb Z(\um)\}$. This is {\em the standard basis} of $I(\lambda)$.
\begin{Prop}\label{radicalofI} For every $\lambda\in\GG$ the $D$-module $I(\lambda)$ has a unique maximal submodule.
\end{Prop}
\pf Let $\Pi$ be the hyperplane of $I(\lambda)$ spanned by all $\{x^{\ui}\eta^{\uj}\otimes 1_{\lambda}\}$ with $|\ui|+|\uj|>0$. We claim that every proper $D$-submodule $M$ of $I(\lambda)$ lies in $\Pi$. If not, then
$$1\otimes 1_{\lambda}+\sum_{|\ui|+|\uj|>0}c_{\ui,\uj}x^{\ui}\eta^{\uj}\otimes 1_{\lambda}\in M\;\text{for some}\;c_{\ui,\uj}\in\k$$
By an argument verbatim to one used in the proof of Lemma \ref{algebraA} the $x_i,\eta_j$ generate a nilpotent subalgebra of $D$. Therefore\newline $z=\sum_{|\ui|+|\uj|>0}c_{\ui,\uj}x^{\ui}\eta^{\uj}\otimes 1_{\lambda}$ is nilpotent, hence $1\otimes 1_{\lambda}\in M$,\newline a contradiction.\qed

Let $R(\lambda)$ denote the radical of $I(\lambda)$, i.e the sum of all proper submodules or zero, if $I(\lambda)$ is simple. We set
$$L(\lambda)=I(\lambda)/R(\lambda).$$
Let $M$ be  a $D$-module. An element $0\neq v\in M$ is called {\em weight element} of weight $\mu\in\GG$ if 
$$g\gamma.v=\mu(g\gamma)v\;\text{holds for all}\;g\gamma\in\Gamma.$$
We say that $v\neq 0$ is {\em primitive} if $J(A).v=0$. We note that every $D$-module $M$ contains primitive weight elements, in fact, a simple submodule of the $A$-socle of $M$ is spanned by a primitive element.

The subalgebra of $D(H)$ generated by $H^*$ and $x_k$ is an Ore extension with the automorphism $\phi_{x_k}:\alpha\mapsto (a_k^{-1}\lh\alpha),\alpha\in H^*,$ and a right $\phi_{x_k}$-derivation $\delta_{x_k}$ determined on generators of $H^*$ by relations \eqref{gammax} and \eqref{etacommx}. Therefore
for every $\alpha\in H^*$, $\alpha x_k^s=\sum_{i=0}^sx_k^i\alpha_i$ for some $\alpha_i\in H^*$. The next two lemmas give a more precise form of these identities.
\begin{Lem}\label{etacommx'}There are polynomials $h_i^{(s)}(t)$ for $s=1,2,\ldots,i=1,2,\ldots s$ such that
\begin{equation}[x_k^s,\eta_k]_{\phi_{x_k}^s}=\sum_{i=1}^sh_i^{(s)}(q_k)x_k^{s-i}\eta_k^{i+1}\end{equation}
\end{Lem}
\pf To simplify notation we drop the subscript $k$ and set $\phi=\phi_{x_k}$. We induct on $s$ noting that the case $s=1$ holds by \eqref{etacommx}. First we have $\phi^s(\eta)=a^{-s}\lh\eta=q^{-s}\eta$ by \eqref{ghiteta}. Now the induction hypothesis and Lemma \ref{xdata}(2) give the identity
$$[x^{s+1},\eta]_{\phi^{s+1}}=q^{-s}x^s[x,\eta]_{\phi}+(\sum_{i=1}^sh_i^{(s)}(q)x^{s-i}\eta^{i+1})x$$
It remains to pass $\eta^{i+1}$ over $x$. We have by Lemma \ref{xdata}(1) that\newline $\eta^{i+1}x=q^{-(i+1)}x\eta^{i+1}+\delta_x(\eta^{i+1})$. There $\delta_x$ is a right $\phi$-derivation with $\delta_x(\eta)=(q-1)\eta^2$ by \eqref{etacommx}. We claim that for every $m$ 
\begin{equation}\label{xdereta}\delta_x(\eta^m)=(q-1)(m)_{q^{-1}}\eta^{m+1}\end{equation}
For by \eqref{prodrule2} $\delta_x(\eta^{m+1})=\delta_x(\eta^m)q^{-1}\eta+\eta^m(q-1)\eta^2$, and assuming the formula for $m$ we get $\delta_x(\eta^{m+1})=(q-1)(q^{-1}(m)_{q^{-1}}+1)\eta^{m+2}$ as asserted.\qed
% polynomials $h_i^{(s)}$ satisfy the following recursive equations. $h_1^{(s+1)}=h_1^{(s)}q^{-2}+q^{-s}(q-1), h_i^{(s+1)}=h_{i-1}^{(s)}(q-1)(i)_{q^{-1}}+h_i^{(s)}q^{-(i+1)}$ and $h_{s+1}^{(s+1)}=h_s^{(s)}(q-1)(s+1)_{q^{-1}}$.
\begin{Lem}\label{gammacommx'}There are polynomials $g_i^{(s)}(t)$ for $s=1,2,\ldots,i=1,2,\ldots s$ such that
\begin{equation}[x_k^s,\gamma]_{\phi_{x_k}^s}=\gamma(a_k^{-s})(x_k^s+\sum_{i=1}^sg_i^{(s)}(q_k)c_k^i(\gamma)x_k^{s-i}\eta_k^i)\gamma\end{equation}
where $c_k^i(\gamma)$ are functions defined in \eqref{functionsofgamma}.
\end{Lem}
\pf We suppress the subscript $k$ as in the preceding lemma. We argue by induction on $s$ starting with \eqref{gammax}. Since\newline $\phi^s(\gamma):=a^{-s}\lh\gamma=\gamma(a^{-s})\gamma$, Lemma \ref{xdata} gives
\begin{align*}[x^{s+1},\gamma]_{\phi^{s+1}}&=x^s\gamma(a^{-s})[x,\gamma]_{\phi}+[x^s,\gamma]_{\phi^s}x\\
&=\gamma(a^{-(s+1)})qc^1(\gamma)x^s\eta\gamma+\gamma(a^{-s})(\sum_{i=1}^sg^i(q)c^i(\gamma)x^{s-i}\eta^i\gamma)x\end{align*}
the last equality by \eqref{gammax} and the induction hypothesis. The proof will be completed if we show that $$(\eta^i\gamma)x=\gamma(a^{-1})(\kappa_1x\eta^i\gamma+\kappa_2(q^i\gamma(ab)-1)\eta^{i+1}\gamma),$$
where $\kappa_i$ are some polynomials of $q$. This equality is derived as follows. First, Lemma \ref{xdata} lets us write\newline $ (\eta^i\gamma)x=x(a^{-1}\lh\eta^i\gamma)+\delta_x(\eta^i\gamma)$. Next, the `$\lh$' action is an algebra homomorphism, hence $a^{-1}\lh\eta^i\gamma=q^{-i}\gamma(a^{-1})\eta^i\gamma$ with the help of Lemma \ref{gconjug}. Lastly, recalling that $\delta_x$ is a right $\phi$-derivation we compute
\begin{align*}\delta_x(\eta^i\gamma)&=\delta_x(\eta^i)\phi(\gamma)+\eta^i\delta_x(\gamma)\\
&=(q-1)(i)_{q^{-1}}\eta^{i+1}\gamma(a^{-1})\gamma+\gamma(a^{-1})qc^1\eta^{i+1}\gamma\\
&=\gamma(a^{-1})[(q-1)(i)_{q^{-1}}+q(\gamma(ab)-1)]\eta^{i+1}\gamma\end{align*}
where the second line is written by \eqref{xdereta} and the basic relation \eqref{gammax}. It remains to note that the expression in the square brackets equals $q^{-(i-1)}(q^i\gamma(ab)-1)$.\qed

We record for the future reference
\begin{equation}\label{xderetagamma}\delta_{x_k}(\eta_k^i\gamma)=q_k^{-(i-1)}\gamma(a_k^{-1})(q_k^i\gamma(a_kb_k)-1)\eta_k^{i+1}\gamma\end{equation}

We move on to the general case of the preceding lemma. For an $n$-tuple $\ui$ we write $c_{\ui}=c_1^{i_1}\cdots c_n^{i_n}$. We recall that $a^{\us}$ stands for $a_1^{s_1}\cdots a_n^{s_n}$.
\begin{Lem}\label{gammaaction} For every pair $(\us,\ut)$ with $\us,\ut\in\mathbb Z(\um)$ there are polynomials $g_{\ui}=g_{\ui}(q_1,\ldots,q_n),\ui\le\us$ such that
\begin{equation}\gamma(x^{\us}\eta^{\ut})\otimes 1_{\lambda}=\lambda(\gamma)\gamma(a^{-\us}b^{\ut})(x^{\us}\eta^{\ut}\otimes 1_{\lambda}+\sum_{\ui\le\us} g_{\ui}c_{\ui}(\gamma)x^{\us-\ui}\eta^{\ut+\ui}\otimes 1_{\lambda})\end{equation}
\end{Lem}
\pf We derive from the previous lemma that $\gamma(x^{\us}\eta^{\ut})\otimes 1_{\lambda}$ is the sum of monomials $m_{\ui}=\gamma(a^{-\us})g_1^{i_1}\cdots g_n^{i_n}c_{\ui}(\gamma)x_1^{s_1-i_1}\eta_1^{i_1}\cdots x_n^{s_n-i_n}\eta_n^{i_n}\gamma\eta^{\ut}$. Now observe that $\eta_i$ skew commutes with every $x_j ,\,j\neq i$, all $\eta_j$, and $\gamma\eta^{\ut}=\gamma(b^{\ut})\eta^{\ut}$ by  Proposition \ref{multipH*'}(1). Thus $m_{\ui}$ can be rewritten as $q_1^{p_1}\cdots q_n^{p_n}x^{\us-\ui}\eta^{\ut+\ui}$ for a suitable integers $p_i$ and the lemma follows.\qed

We begin to build the weight space decomposition of $I(\lambda)$. For $g\in G$ we denote by $\hat g$ the character of $\G$ sending $\gamma$ to $\gamma(g)$ for every $\gamma\in\G$. For $\us,\ut\in(\mathbb Z)^n$ and $\lambda\in\Gamma$ we define the character $\lambda_{\us,\ut}$ by 
$$\lambda_{\us,\ut}=\lambda\widehat{a^{-\us}b^{\ut}}\chi^{\us+\ut}.$$
Recall the idempotent $e_{\lambda}=|\Gamma|^{-1}\sum_{g\gamma}\lambda^{-1}(g\gamma)g\gamma$ associated to $\lambda\in\Gamma$. We denote by $e_{\us,\ut}$ the idempotent corresponding to $\lambda_{\us,\ut}$. We define vector $v_{\us,\ut}$ by the formula $v_{\us,\ut}=e_{\us,\ut}(x^{\us}\eta^{\ut}\otimes 1_{\lambda})$. In the same spirit we let $I_{\us,\ut}(\lambda)$ denote the weight space $e_{\us,\ut}I(\lambda)$. We put an equivalence relation on the set $\mathbb Z(\um)\times\mathbb Z(\um)$ by declaring  $(\us,\ut)\sim (\us',\ut')\;\text{if and only if}\;\lambda_{\us,\ut}=\lambda_{\us',\ut'}$. We let $[\us,\ut]$ denote the equivalence class of $(\us,\ut)$.
\begin{Lem}\label{weightdecomp}
\begin{itemize}
\item[(1)] The set 
$$\{v_{\us,\ut}|\us,\ut\in\mathbb Z(\um)\}$$
\;\text{is a basis for}\; $I(\lambda)$.

\item[(2)] The set 
$$\{v_{\us',\ut'}|(\us',\ut')\in[\us,\ut]\}$$
is a basis for $I_{\us,\ut}(\lambda)$.
\item[(3)] $I(\lambda)=\bigoplus \{I_{\us,\ut}|\us,\ut\in\mathbb Z(\um)\}$
\end{itemize}
\end{Lem}
\pf (1) The defining relations of $H$ and \eqref{geta} make it clear that\newline $x^{\us}\eta^{\ut}\otimes 1_{\lambda}$ has $G$-weight $\lambda|_G\chi^{\us+\ut}$. This observation combined with the previous lemma shows that 
\begin{equation}\label{weightvector}v_{\us,\ut}=x^{\us}\eta^{\ut}\otimes 1_{\lambda}+\sum_{\ui<\us}g_{\ui}\overline{c_{\ui}}x^{\us-\ui}\eta^{\ut+\ui}\otimes 1_{\lambda}\end{equation}
where $\overline{c_{\ui}}:=|G|^{-1}\sum_{\gamma\in\G}c_{\ui}(\gamma)$ is the average value of $c_{\ui}$ over $\Gamma$. We see that $v_{\us,\ut}$ has the leading term $x^{\us}\eta^{\ut}\otimes 1_{\lambda}$. We claim that the set in part (1) is linearly independent. Else, there is a linear relation\newline (*) $\sum\kappa_{\us,\ut}v_{\us,\ut}=0,\;0\neq\kappa_{\us,\ut}\in\k$. Pick $(\us',\ut')$ such that $\us'$ is the largest among all $(\us,\ut)$ involved in (*) and $\ut'$ is the smallest among all $(\us',\ut)$ involved in (*) in the ordering of \S 2.1. Then $x^{\us'}\eta^{\ut'}\otimes 1_{\lambda}$ can not get cancelled, a contradiction. As $\dim\,I(\lambda)=\text{card}\,|\mathbb Z(\um)|$, the assertion holds.

(2) and (3). We note that every $v_{\us',\ut'}$ with $(\us',\ut')\in[\us,\ut]$ lies in $I_{\us,\ut}(\lambda)$ by the definition of the latter, hence $\dim\,I_{\us,\ut}(\lambda)\ge |[\us,\ut]|$. However, the sum of cardinalities of the distinct sets $[\us,\ut]$ equals $\dim\,I(\lambda)$. On the other hand, sum of $I_{\us,\ut}(\lambda)$ is direct of dimension no greater than $\dim\,I(\lambda)$. This proves (2) and (3).\qed

The next theorem is the main result of this section. We note that the theorem and its proof resemble theorems of Curtis \cite{Cu} and Lusztig \cite{Lu2} on parametrization of simple modules.
\begin{Thm}\label{simpleD-mod} The simple modules $L(\lambda),\,\lambda\in\widehat\Gamma$, are a full set of representatives of simple $D$-modules.
\end{Thm}
\pf As every simple $D$-module $L$ is generated by a primitive weight element, $L$ is the image of $I(\lambda)$ for a suitable $\lambda$. It remains to show that $L(\lambda)\simeq L(\mu)$ implies $\lambda=\mu$. Let $\Pi$ be the hyperplane in $I(\lambda)$ of Proposition \ref{radicalofI}, and denote by $P$ its image in $L(\lambda)$. Were $L(\lambda)\simeq L(\mu)$, there would be a primitive vector $v$ of weight $\mu$ in $L(\lambda)$.

We note that $P$ is a proper subspace of $L(\lambda)$. For, a linear relation
$$1\otimes 1_{\lambda}+\sum c_{\ui,\uj}x^{\ui}\eta^{\uj}\otimes 1_{\lambda}\equiv 0\mod R(\lambda)$$
implies $1\otimes 1_{\lambda}\in R(\lambda)$, because the $x_l,\eta_k$ generate a nilpotent subalgebra, a contradiction. As $\lambda\neq\mu$ by assumption, $v$ is a linear combination of the images of $v_{\us,\ut}$ with $(\us,\ut)\neq (\uo,\uo)$. By \eqref{weightvector} $v\in P$, hence, as $v$ is primitive, $L(\lambda)$ is the span of all $x^{\ui}\eta^{\uj}v$. Refering to Lemma \ref{etacommx'} we see that $P$ is invariant under multiplication by $\eta_k$ and, of course, $x_l$. Thus $L(\lambda)\subset P$, a contradiction.\qed
\pagebreak
%%%%%%%%%%%%%%%%%%%%%%%%%%%%%%%%%%%%%%%%%%%%%%%%%%%%%%%%%%%%%%%%%%%%%%%%%%%%%%%%%%%%%%%%%%%%%%%%%%%%%%%%%%%%%%%%%%%%%%%%%%%%%%%%%%%%%%%%%%%%%%%%%%%%%%%%

\subsection{The Loewy Filtration of $I(\lambda)$}\label{Lfiltration}
\subsubsection{Action of generators on standard basis}
Our ultimate goal is to describe action of generators of $D$ on the weight basis of $I(\lambda)$. The next proposition is a key step toward this goal. It runs smoothly under a certain restriction on the datum for $H$. We don't know whether this restriction is essential for our results.

For an even $m$ we put $m'=m/2$. 
\begin{Def}\label{halfcleandata} We say that a simply linked datum $\mathcal D$ is {\em half-clean} if $G$ does not have nontrivial relations $r=1$ of the form
$$r=\prod_{i=1}^n(a_ib_i)^{t_i}$$
with $t_i=0\;\text{or}\; m_i'$.

The above definition has the following implication for the set of relations of $G$. $G$ does not have relators $r=\prod_{i=1}^n(a_ib_i)^{t_i}\neq 1\;\text{with}\;t_i< m_i$ for all $i$. For, since $\chi_j(a_ib_i)=1$ for every $j\neq i$, $\chi_i(r)=\chi_i(a_ib_i)^{t_i}=q_i^{2t_i}$. Therefore $r=1$ implies $2t_i\equiv 0\mod m_i$. Thus $t_i=0$ if $m_i$ is odd, and $t_i=0\;\text{or}\; m'_i$, otherwise. 
\end{Def}
\begin{Prop}\label{e-action} Suppose $\mathcal D$ is a half-clean datum, $\us\in\mathbb Z(\um)$ and $\ut\in{\mathbb Z}^n$. 
\begin{itemize}
\item[(1)] For every $\ui,\,\uo<\ui\le\us$, $$e_{\us,\ut}(x^{\us-\ui}\eta^{\ut+\ui}\otimes 1_{\lambda})=0.$$
\item[(2)] For every $k\in\un$ and every $\ui,\,\uo\le\ui\le\us$, $$e_{\us+u_k,\ut}(x^{\us-\ui}\eta^{\ut+u_k+\ui}\otimes 1_{\lambda})=0.$$
\item[(3)] For every $i,0\le i\le s_k-1$, $$e_{\us,\ut-u_k}(x^{\us-(i+1)u_k}\eta^{\ut+iu_k}\otimes 1_{\lambda})=0.$$
\end{itemize}
\end{Prop}
\pf As a preliminary to the proof we point out that for every abelian group $G$, $\sum_{\gamma\in\G}\gamma(g)=|G|\delta_{1,g}$.

(1) Put $v=x^{\us-\ui}\eta^{\ut+\ui}\otimes 1_{\lambda}$. By Lemma \ref{gammaaction} and \eqref{geta} $v$ acquires weight $\lambda_{\us-\ui,\ut+\ui}$ upon multiplication by $g\gamma$. One can see readily that $\lambda^{-1}_{\us,\ut}\lambda_{\us-\ui,\ut+\ui}=\widehat{(ab)^{\ui}}$. Consequently Lemma \ref{gammaaction} gives the equality
\begin{equation}e_{\us,\ut}.v=\sum_{\uj\le\us-\ui} g_{\uj}\{|G|^{-1}\sum_{\gamma\in\G}\gamma((ab)^{\ui})c_{\uj}(\gamma)\}x^{\us-\ui-\uj}\eta^{\ut+\ui+\uj}\otimes 1_{\lambda}\tag{*}\end{equation}
We begin to work out the inner sum in (*). First, for a positive integer $m$ and every $k\in\in$ we have the expansion
$$c_k^m(\gamma)=\sum_{l=0}^m(-1)^{m-l}\binom ml_{q_k}q_k^{\binom l2}\gamma((a_kb_k)^l)$$
by the Gauss' binomial formula \cite{KC}.

Replacing $m$ by $j_k$ and $l$ by $l_k$ and letting $k$ run over $\un$ we see that the inner sum in (*) equals $|G|^{-1}\sum_{\ul\le\uj}\gamma((ab)^{\ui+\ul})$. Since $\ui+\ul\le\us<\um$ our assumption on $\mathcal D$ imply that $(ab)^{\ui+\ul}\neq 1$ for every $\ul$, whence the sum vanishes by the opening remark, and this proves (1).

(2) We set $v=x^{\us-\ui}\eta^{\ut+u_k+\ui}\otimes 1_{\lambda}$. A simple verification gives\newline $\lambda^{-1}_{\us+u_k,\ut}\lambda_{\us-\ui,\ut+u_k+\ui}=\widehat{(ab)^{\ui+u_k}}$. Repeating the argument leading up to the equality (*) we derive that 
$$e_{\us+u_k,\ut}v=\sum_{\uj\le\us-\ui}g_{\uj}\{|G|^{-1}\sum_{\gamma\in\G}\gamma((ab)^{\ui+u_k})c_{\uj}(\gamma)\}x^{\us-\ui-\uj}\eta^{\ut+u_k+\ui+\uj}\otimes 1_{\lambda}.$$
Applying the Gauss' binomial formula to $c_{\uj}(\gamma)$ we have
\begin{align}\sum_{\gamma\in\G}\gamma((ab)^{\ui+u_k})c_{\uj}&=\sum_{\ul\le\uj}\kappa_{\ul}\sigma_{\ul},\kappa_{\ul}\in\k,\tag{**}\quad\text{where}\\
\sigma_{\ul}&=\sum_{\gamma\in\G}\gamma((ab)^{\ui+u_k+\ul})\tag{***}\end{align}
Further we note that, as $\ui+\ul\le\us$, if $s_k<m_k-1$, we have $\uo<\ui+u_k+\ul<\um$, hence the sum (***) is zero for all $\ul$, and therefore the sum (**) is zero for all $\uj$. Same conclusion holds if $s_k=m_k-1$, but $j_k<s_k$. What remains to be considered is the case when $l_k=j_k=m_k-1-i_k$. Now $l_k+i_k+1=m_k$, hence the exponent on $\eta_k$ in $\eta^{\ut+u_k+\ui+\uj}$ is $\ge m_k$, and therefore this term vanishes. Thus the sum (**) is always zero, and (2) is done.

(3) can be reduced to (1) by replacing $\ut$ with $\ut-u_k$.\qed

From now on we assume that $\mathcal D$ is half-clean. We need one more commutation formula.
\begin{Lem}\label{xicommx'}For every $s\ge 1$ and every $k\in\un$ there are functions $r_i^s$ of $q_k,\;1\le i\le s$ such that
\begin{equation}[x_k^s,\xi_k]_{\phi_{x_k}^s}=(s)_{q_k}x_k^{s-1}(a_k-q_k^{-(s-1)}\chi_k)+\sum_{i=1}^{s-1}r_i^sx_k^{s-1-i}\eta_k^i\chi_k\end{equation}
\end{Lem}
\pf As before we drop the subscript $k$ throughout.  We induct on $s$ noting that the formula holds for $s=1$ by \eqref{xicommx}. To carry out the induction step we use Lemma \ref{xdata}(2), namely
\begin{align*}[x^{s+1},\xi]_{\phi^{s+1}}&=x^s[x,a^{-s}\lh\xi]_{\phi}+[x^s,\xi]_{\phi^s}x\\
\intertext{(which by \eqref{ghitxi},\eqref{xicommx}, and the induction hypothesis equals)}
&x^s(a-\chi)+(s)_qx^{s-1}(a-q^{-(s-1)}\chi)x+\sum_{i=1}^{s-1}r_i^sx^{s-1-i}\eta^i\chi x\end{align*}
For every $m=s-i,i>0$ the coefficient of $x^m$ is of the form $r_i^{s+1}\eta^{i+1}\chi$ by \eqref{xderetagamma} as needed. It remains to compute the coefficient of $x^s$. Using \eqref{gammax}, we have
$\chi x=q^{-1}x\chi+(q^2-1)\eta\chi$, hence $(a-q^{_(s-1)}\chi)x=x(qa-q^{-s}\chi)+\kappa\eta\chi$ for some $\kappa\in\k$. Thus this coefficient equals $(a-\chi)+(s)_q(qa-q^{-s}\chi)=(s+1)_q(a-q^{-s}\chi)$ and the proof is complete.\qed

In what follows we set $v_{\us,\ut}=0$ if $s_k$ or $t_k$ is zero for some $k\in\un$. We derive now action of $\eta_k,\xi_k$ on $I(\lambda)$.
\begin{Prop}\label{eta&xiaction} For each $k\in\un$ and $\us,\ut\in\mathbb Z(\um)$ there are roots of unity $\theta,\theta'$ such that 
\begin{itemize}
\item[(1)] $\eta_kv_{\us,\ut}=\theta v_{\us,\ut+u_k}$
\item[(2)] $\xi_kv_{\us,\ut}=\theta'(s_k)_{q_k}(\lambda(a_k\chi_k^{-1})-q_k^{-(s_k-1)})v_{\us-u_k,\ut}$
\end{itemize}
\end{Prop}
\pf (1) We start off by deriving action of $\eta_k,\xi_k$ on idempotents $e_{\us,\ut}$. Pick $\mu\in\Gamma$. We have
\begin{align*}\eta_ke_{\mu}=|\Gamma|^{-1}\sum\mu((g\gamma)^{-1})\eta_kg\gamma&=|\Gamma|^{-1}(\sum\mu((g\gamma)^{-1})\chi_k(g^{-1})\widehat b_k(\gamma^{-1})\\&=e_{\nu}\eta_k\quad\text{where}\;\nu=\mu\widehat b_k\chi_k\end{align*}
Setting $\mu=\lambda_{\us,\ut}$ we obtain the first identity, namely
\begin{equation}\label{etaidempotent}\eta_ke_{\us,\ut}=e_{\us,\ut+u_k}\eta_k\end{equation}
For $\xi_k$ we can show similarly that
\begin{equation}\label{xiidempotent} \xi_ke_{\us,\ut}=e_{\us-u_k,\ut}\xi_k\end{equation}
We continue with part (1). In order to apply \eqref{etaidempotent} we must expand $\eta_k x^{\us}\eta^{\ut}\otimes 1_{\lambda}$ in the standard basis of $I(\lambda)$. By \eqref{etacommx} $\eta_k$ skew commute with $x_l$ for all $l\neq k$. For $l=k$ we use Lemma \ref{etacommx'} to pass $\eta_k$ over $x_k^{s_k}$. Noting that $\eta_k$ skew commutes with every $\eta_l$ we arrive at the equality
\begin{equation}\label{etamonomial}\eta_kx^{\us}\eta^{\ut}\otimes 1_{\lambda}=\theta x^{\us}\eta^{\ut+u_k}\otimes 1_{\lambda}+\sum_{i=1}^sh_ix^{\us-iu_k}\eta^{\ut+(i+1)u_k}\otimes 1_{\lambda}\end{equation} 
where $h_i$ are some elements of $\k$. By part (1) of the preceding proposition $e_{\us,\ut+u_k}$ annihilates $x^{\us-iu_k}\eta^{\ut+(i+1)u_k}\otimes 1_{\lambda}$ for every $i\ge 1$, which gives (1).

(2) We begin by working out an expansion of $\xi_kx^{\us}\eta^{\ut}\otimes 1_{\lambda}$ in the standard basis of $I(\lambda)$. Set $\us'=(s_1,\ldots,s_{k-1},0,\ldots,0)$ and $\us''=(0,\ldots,0,s_{k+1},\ldots,s_n)$. Noting that $\xi_k$ commutes with $x_l$ for every $l\neq k$ and commutes with all $\eta_l$, as well as the equality $\xi_k.1_{\lambda}=0$, we have with the help of Lemma \ref{xicommx'} the first equality,viz.
\begin{equation}\xi_kx^{\us}\eta^{\ut}\otimes 1_{\lambda}=(s_k)_{q_k}w_0+\sum_{i=1}^{s_k-1}\kappa_iw_i\tag{*}\end{equation}
where $\kappa_i\in\k$, $w_0=x^{\us'}x_k^{s_k-1}(a_k-q_k^{-(s_k-1)}\chi_k)x^{\us''}\eta^{\ut}\otimes 1_{\lambda}$ and\newline
$w_i=x^{\us'}x_k^{s_k-i-1}\eta_k^i\chi_kx^{\us''}\eta^{\ut}\otimes 1_{\lambda}$.

The rest of the proof will be carried out in steps. (i) By \eqref{orthogonality} and \eqref{gammax} $\chi_kx_l=q_{lk}^{-1}x_l\chi_k$, and also $a_kx_l=q_{lk}^{-1}x_la_k$ for $l\neq k$. Further $a_k$ skew commutes with $\eta^{\ut}$ with scalar $\chi^{\ut}(a_k)$ and $\chi_k$ skew commutes with $\eta^{\ut}$ with scalar $\chi_k(b^{\ut})$. As $\chi^{\ut}(a_k)=\chi_k(b^{\ut})$ by condition (D1), it follows that $w_0=\theta'(\lambda(a_k\chi_k^{-1}-q^{-(s_k-1)})x^{\us-u_k}\eta^{\ut}\otimes 1_{\lambda}$.

(ii) We claim that $\eta_k^i\chi_k$ skew commutes with $x_l$ for every $l\neq k$, or, equivalently, $\delta_{x_l}(\eta_k^i\chi_k)=0$. Indeed, $\delta_{x_l}(\eta_k)=0$ holds by \eqref{etacommx} and $\delta_{x_l}(\chi_k)=0$ follows from the first line of proof of (i). Since $\delta_{x_l}$ is a skew derivation, the claim follows. As $\eta_k,\chi_k$ skew commute with all $\eta_l$ it becomes clear that $w_i=\kappa x^{\us-(i+1)u_k}\eta^{\ut+iu_k}\otimes 1_{\lambda}$ for some $\kappa\in\k$. 

From (i) and (ii) we have the expansion 
\begin{align}\label{ximonomial}\xi_kx^{\us}\eta^{\ut}\otimes 1_{\lambda}&=\theta'(s_k)_{q_k}\lambda(a_k\chi_k^{-1}-q_k^{-(s_k-1)})x^{\us-u_k}\eta^{\ut}\otimes 1_{\lambda}\nonumber\\&+\sum_{i=1}^{s_k-1}\kappa_ix^{\us-(i+1)u_k}\eta^{\ut+iu_k}\otimes 1_{\lambda}\end{align}

Applying $e_{\us-u_k,\ut}$ to equality \eqref{ximonomial}  we arrive by Proposition \ref{e-action}(1) at the desired result.\qed

The next lemma gives a commutation relation between some generators of $H$ and primitive idempotents.
\begin{Lem}\label{xyidempotents} For every $\in\un$ and $\us,\ut\in\mathbb Z(\um)$
\begin{itemize}\item[(1)] $x_ke_{\us,\ut}=e_{\us+u_k,\ut}x_k+q_k(e_{\us,\ut+u_k}-e_{\us+u_k,\ut})\eta_k$
\item[(2)] $y_ke_{\us,\ut}=e_{\us,\ut-u_k}y_k+q_k(e_{\us,\ut-u_k}-e_{\us-u_k,\ut})b_k\xi_k$
\end{itemize}
\end{Lem}
\pf  Combining Lemmas \ref{xdata}, \ref{gconjug}, \ref{gammacommutators} and Proposition \ref{multipH*'} we compute
\begin{align*} x_kg\gamma&=\chi_k(g^{-1})\gamma(a_k)g\gamma x_k -q_k\chi_k(g^{-1})(\gamma(a_k)-\gamma(b_k^{-1}))g\gamma\eta_k\\
y_kg\gamma&=\chi_k(g)\gamma(b_k)g\gamma y_k+q_k\chi_k(g)(\gamma(b_k)-\gamma(a_k^{-1}))g\gamma b_k\xi_k\end{align*}
It follows that for every $\mu\in\GG$ $x_ke_{\mu}$ splits up into the sum\newline $x_ke_{\mu}=e_{\mu'}x_k-q_k(e_{\mu'}-e_{\mu''})\eta_k$ where $\mu'=\mu\chi_k\widehat{a_k^{-1}}$ and $\mu''=\mu\chi_k\widehat{b_k}$. Similarly $y_ke_{\mu}=e_{\mu'}y_k+q_k(e_{\mu'}-e_{\mu''})b_k\xi_k$ where $\mu'=\mu\chi_k^{-1}\widehat{b_k^{-1}}$ and $\mu''=\mu\chi_k^{-1}\widehat{a_k}$. Setting $\mu=\lambda_{\us,\ut}$ we obtain the lemma.\qed

We can now describe action of $x_k$ on $I(\lambda)$. Briefly, $x_k$ acts as a {\em raising} operator, however, not literally, because the set of weights is not ordered.
\begin{Prop}\label{x-action} For every $k\in\un$ and $\us,\ut\in\mathbb Z(\um)$ there are roots of unity $\theta,\theta'$ such that
$$x_kv_{\us,\ut}=\theta v_{\us+u_k,\ut}+\theta'v_{\us,\ut+u_k}.$$
\end{Prop}
\pf Since $x_k$ skew commutes with every $x_l,l\neq k,\,x_k(x^{\us}\eta^{\ut})\otimes 1_{\lambda}=\theta x^{\us+u_k}\eta^{\ut}\otimes 1_{\lambda}$. By the previous lemma the latter monomial contributes $v_{\us+u_k,\ut}$ to $x_kv_{\us,\ut}$. We move on to $w:=\eta_kx^{\us}\eta^{\ut}\otimes 1_{\lambda}$. From the expansion \eqref{etamonomial} and Proposition \ref{e-action}(1) we see that $e_{\us,\ut+u_k}w=\theta'v_{\us,\ut+u_k}$. On the other hand part (2) of that proposition gives $e_{\us+u_k,\ut}w=0$, and this completes the proof.\qed

We isolate one step of the next proposition in-
\begin{Lem}\label{ypowereta} For every $k\in\un$ and $t<m_k$ there holds
$$y_k\eta_k^t=q_k^t\eta_k^ty_k-(t)_{q_k}\eta^{t-1}(q_k^t\chi_kb_k-\epsilon)$$
\end{Lem}
\pf We use induction on $t$, the case $t=1$ is covered by Lemma \ref{etacommutators}. Let $\phi:H^*\to H^*$ be the automorphism $\phi(\alpha)=b_k\lh\alpha,\alpha\in H^*$. Below we drop the subscript $k$ on $y_k,\eta_k$. For the induction step we make use of \eqref{prodrule1} which allows us to write
\begin{align*}{}_{\phi}[y,\eta^{t+1}]&={}_{\phi}[y,\eta^t]\eta+q^t\eta^t{}_{\phi}[y,\eta]\\
&=-(t)_q\eta^{t-1}((q^t\chi b-\epsilon)\eta-q^t\eta^t(\chi b-\epsilon))\end{align*}
where for the last line we used the induction hypothesis and the basis of induction. One can check easily the equality $\chi b\eta=q^2\eta\chi b$ which leads to the equality  ${}_{\phi}[y,\eta^{t+1}]=-\eta^t(q^t((t)_qq+1)\chi b-((t)_q+q^t)\epsilon)=-\eta^t(t+1)_q(q^t\chi b-\epsilon)$.\qed

We finish this section by showing that $y_k$ acts as a lowering operator. More presicely
\begin{Prop}\label{y-action}For every $k\in\un$ and $\us,\ut\in\mathbb Z(\um)$ there are roots of unity $\theta,\theta'$ such that
\begin{align}\label{yaction}y_kv_{\us,\ut}&=\theta(t_k)_{q_k}(\lambda(\chi_kb_k)-q_k^{-(t_k-1)})v_{\us,\ut-u_k}\\&+\theta'(s_k)_{q_k}(\lambda(a_k\chi_k^{-1})-q_k^{-(s_k-1)})v_{\us-u_k,\ut}\nonumber\end{align}
\end{Prop}
\pf In view of Lemma \ref{xyidempotents} we want to express $y_kx^{\us}\eta^{\ut}\otimes 1_{\lambda}$ and $\xi_kx^{\us}\eta^{\ut}\otimes 1_{\lambda}$ in the standard basis of $I(\lambda)$. For the second of those vectors the expansion is given by \eqref{ximonomial}. Multiplying \eqref{ximonomial}
by $e_{\us-u_k,\ut}$ and $e_{\us,\ut-u_k}$ in turn, we get the second summand of \eqref{yaction} and zero, respectively, by Proposition \ref{e-action}(1),(3).

We turn now to $w=y_kx^{\us}\eta^{\ut}\otimes 1_{\lambda}$. $y_k$ skew commutes with $x_l$ and $\eta_l$ for $l\neq k$. For $l=k$ the product $y_kx_k^{s_k}$ can be streightened out by Lemma \ref{commutatorrel1}(2). In addition $a_kb_k$ commutes with $x_l$ for $l\neq k$ and skew commutes with all $\eta_l$. It follows that  
$$ w=\kappa x^{\us}y_k\eta^{\ut}\otimes 1_{\lambda}+\mu x^{\us-u_k}\eta^{\ut}\otimes 1_{\lambda}$$
where $\kappa,\mu\in\k$ with $\kappa$ a root of unity. Since $e_{\us,\ut-u_k}$ annihilates the second monomial in $w$ by Proposition \ref{e-action}(3), we turn to the first summand in $w$. Using the preceding lemma and the fact that $y_k.1_{\lambda}=0$ $w$ reduces to the form $\theta(t_k)_{q_k}(\lambda(\chi_k b_k)-q_k^{-(t_k-1)})x^{\us}\eta^{\ut-u_k}\otimes 1_{\lambda}$ which, upon multiplication by $e_{\us,\ut-u_k}$,  becomes the first summand of \eqref{yaction}.\qed
%%%%%%%%%%%%%%%%%%%%%%%%%%%%%%%%%%%%%%%%%%%%%%%%%%%%%%%%%%%%%%%%%%%%%%%%%%%%%%%%%%%%%%%%%%%%%%%%%%%%%%%%%%%%%%%%%%%%%%%%%%%%%%%%%%%%%%%%%%%%%%%%%%%%%%

\subsubsection{The Main Theorems}\label{maintheorems}
In this section we will show that representation theory of $D(H)$ in $I(\lambda)$ follows the pattern established for $H$ in Theorem \ref{Lseries1}. 
\begin{Def}\label{singularvalues} For $\lambda\in\GG$ we let $S(\lambda)$ stand for all $j\in\un$ satisfying the condition
\begin{align} \lambda(a_j\chi_j^{-1})&=q_j^{-e_j}\label{sing1}\;\text{or}\\\lambda(b_j\chi_j)&=q_j^{-e'_j}\label{sing2}\end{align}
for some $0\le e_j,e_j'\le m_j-2$
\end{Def}
Those constants depend on $\lambda$. We write them as $e_j(\lambda),e'_j(\lambda)$ when this dependance must be emphasized.

For an $n$-tuple $\ua\in\mathbb Z^n$ we define {\em the rank} of $\ua$ as the number $\text{rk}(\ua)$ of all $j\in S(\lambda)$ satisfying $a_j\ge e_j+1$. We define {\em the rank} of $x^{\us}\eta^{\ut}\otimes 1_{\lambda}$ by $\text{rk}(x^{\us}\eta^{\ut}\otimes 1_{\lambda})=\text{rk}(\us)+\text{rk}(\ut)$.
\begin{Lem}\label{primitive} A weight vector $v_{\us,\ut}$ is primitive if and only if\newline $s_j=0,e_j+1$ or $t_j=0,e'_j+1$ for every $j\in S(\lambda)$, and $s_k=0,t_k=0$ for every $k\notin S(\lambda)$.
\end{Lem}
\pf This follows immediately from Proposition \ref{eta&xiaction} and Proposition \ref{y-action}.\qed

Let's denote the $D$- module generated by a set $X$ by $\Dg X$.
\begin{Prop}\label{radicalI} $R(\lambda)$ is generated by the primitive elements of rank~1.
\end {Prop}
\pf Suppose $v$ is primitive. From definition of primitivity we have that $\Dg v$ is the span of all $x^{\us}\eta^{\ut}v$. Let now $v$ has rank $1$. Then $v$ has the form $w_j=v_{(e_j+1)u_j,\uo}$ or $w'_j=v_{\uo,(e'_j+1)u_j}$. From Propositions \ref{eta&xiaction} and \ref{x-action} we see that $\Dg v$ is the span of either all $v_{\uc,\ud}$ with $\uc\ge (e_j+1)u_j$ and $\ud\ge\uo$ or $v_{\uc,\ud}$ with $\uc\ge\uo$ and $\ud\ge (e_j'+1)u_j$. Therefore were $R(\lambda)\neq \sum_{j\in S(\lambda)}\Dg {w_j}+\Dg {w_j'}$ there would be a weight vector $v\in R(\lambda)$ involving $v_{\us,\ut}$ with $s_j\le e_j$ and $t_j\le e_j'$ for all $j\in S(\lambda)$. We may assume that $v_{\us,\ut}$ is the minimal such in the sense that $R(\lambda)$ doesn't contain weight vectors whose expansion in the standard basis involves $v_{\uc,\ud}$ with either $\uc< \us\;\text{and}\;\ud\le\ut$ or $\uc\le\us\;\text{and}\;\ud< \ut$. Let 
$$v=v_{\us,\ut}+\sum c_{\us',\ut'}v_{\us',\ut'},\;c_{\us',\ut'}\in\k$$
be the expansion of $v$ in the standard basis of $I(\lambda)$. If $\us=\uo=\ut$, then $1\otimes 1_{\lambda}+\sum_{(\us,\ut)\neq (\uo,\uo)}v_{\us,\ut}\in I(\lambda)$, hence, as $x_i,\eta_i$ generate a nilpotent subalgebra, $1\otimes 1_{\lambda}\in R(\lambda)$, a contradiction. Else, set $v'=\xi_k v$ and $v''=y_kv$. Then, either $s_k\neq 0$ for some $k$, hence $v'$ involves a weight vector smaller than $v_{\us,\ut}$ by Prposition \ref{eta&xiaction}, or $\us=\uo$, but $t_k\neq 0$, and then the same holds for $v''$ by Proposition \ref{y-action}. In either case we arrive at a contradiction.\qed

For the dimension formula we introduce some subsets of $S(\lambda)$. We let $S^{(1)}(\lambda),S^{(2)}(\lambda)$ and $S^{(3)}(\lambda)$ be defined by the conditions
\begin{align*}S^{(1)}(\lambda)&=\{j\in S(\lambda)|j\;\text{satisfies only}\;\eqref{sing1}\}\\
S^{(2)}(\lambda)&=\{j\in S(\lambda)|j\;\text{satisfies only}\;\eqref{sing2}\}\\
S^{(3)}(\lambda)&=S^{(1)}(\lambda)\cap S^{(2)}(\lambda)\end{align*}
\begin{Cor}\label{dimension} For every $\lambda\in\GG$
$$\dim L(\lambda)=\prod_{j\notin S(\lambda)}\!\!m_j^2\!\!\prod_{j\in S^{(1)}(\lambda)}\!\!(e_j+1)m_j\!\!\prod_{j\in S^{(2)}(\lambda)}\!\!(e'_j+1)m_j\!\!\prod_{j\in S^{(3)}(\lambda)}\!\!(e_j+1)(e'_j+1).$$
\end {Cor}
\pf By definition $L(\lambda)=I(\lambda)/R(\lambda)$. By the above proposition $L(\lambda)$ is the span of all $v_{\us,\ut}$ with $s_j\le e_j$ and $t_j\le e'_j$ for $j\in S^{(1)}(\lambda)$ or $j\in S^{(2)}(\lambda)$, respectively, and arbitrary integers within $[0, m_j-1]$, otherwise.\qed

We refer to section \ref{repsofH} for a discussion of the Loewy and socle series. We denote by $\ell(\lambda)$ the Loewy length of $I(\lambda)$.  
\begin{Thm}\label{Lseries2} {\rm (1)} For every $m<\ell(\lambda)$ $R^m(\lambda)$ is generated by the primitive vectors of rank $m$.

{\rm (2)} The radical and the socle series coincide.

{\rm (3)} The Loewy layers ${\mathcal L}^m:=R^m(\lambda)/R^{m+1}(\lambda)$ are given by the formula
$${\mathcal L}^m\simeq\oplus\{L(\mu)|\mu\,\text{is weight of primitive basis vector of rank}\;m\}.$$

{\rm (4)} $\ell(\lambda)=|S^{(1)}(\lambda)|+|S^{(2)}(\lambda)|+1$
\end{Thm}
\pf (1)-(4). We abbreviate $R(\lambda)$ to $R$. We induct on $m$. The fact that $I(\lambda)$ is generated by $1\otimes 1_{\lambda}$ gives the basis of induction. Suppose (1) holds for $R^m$. The induction step will be carried out in steps.

(i) Let $v_{\us,\ut}$ be a primitive vector of rank $m$, and $\mu=\lambda_{\us,\ut}$ be its weight. We claim that $S^{(i)}(\mu)=S^{(i)}(\lambda)$ for $i=1,2,3$. Namely, $e_j(\mu)=e_j(\lambda)$ if $s_j=0$, and $e_j(\mu)=m_j-e_j(\lambda)-2$, otherwise. Similarly, $e'_j(\mu)=e'_j(\lambda)$ if $t_j=0$, and $e'_j(\mu)=m_j-e'_j(\lambda)-2$, otherwise. This follows from the calculation
\begin{align*} &\chi^{\us+\ut}\widehat{a^{-\us}b^{\ut}}(a_k\chi_k^{-1}=\chi^{\us+\ut}(a_k)\chi_k(a^{\us}b^{-\ut})\\
&=[\prod_{m\neq k}\chi_m^{s_m}(a_k)\chi_m^{t_m}(a_k)]\chi_k^{s_k}(a_k)\chi_k^{t_k}(a_k)[\prod_{m\neq k}\chi_k(a_m^{s_m})\chi_k(b^{-t_m})]\\&\chi_k(a_k^{s_k})\chi_k(b_k^{-t_k})=\prod_{m\neq k}[\chi_m(a_m)\chi_k(a_m)]^{s_m}[\chi_m(a_k)\chi_k(b_m^{-1})]^{t_m}\\&\chi_k^{2s_k}(a_k)
=\chi_k^{2s_k}(a_k)\end{align*}
with the last equality by the datum conditions (D1-D2). Similarly $\chi^{\us+\ut}\widehat{a^{-\us}b^{\ut}}(b_k\chi_k)=\chi_k^{2t_k}(b_k)$. Therefore $\mu(a_k\chi_k^{-1})=\lambda(a_k\chi_k^{-1})\chi_k^{2s_k}(a_k)$ and $\mu(b_k\chi_k)=\lambda(b_k\chi_k)\chi_k^{2t_k}(b_k)$. It follows that if $s_k=0$ or $t_k=0$, then the value of $\mu$ and $\lambda$ at $a_k\chi_k^{-1}$ or $b_k\chi_k$ coincide. Else, $s_k=e_k(\lambda)+1$ or $t_k=e_k'(\lambda)+1$, whence $\mu(a_k\chi_k^{-1})=\lambda(a_k\chi_k^{-1})q_k^{2(e_k(\lambda)+1)}=q_k^{-(m_k-e_k(\lambda)-2)}$, or similarly $\mu(b_k\chi_k)=q_k^{-(m_k-e'_k(\lambda)-2)}$. This proves our claim.

(ii) Since $v_{\us,\ut}$ is primitive of weight $\mu$ there is a $D$-epimorphism $\phi:I(\mu)\to Dv_{\us,\ut}$ determined on the generator by $\phi: 1\otimes 1_{\mu}\mapsto v_{\us,\ut}$. It follows that $R(Dv_{\us,\ut})=\phi(R(\mu))$, hence, by the preceding proposition, we have that $R(Dv_{\us,\ut})=\sum D\phi(w)$ where $w$ runs over the primitive vectors of $I(\mu)$ of rank $1$. Put $f_j=e_j(\mu)+1$ and $f_j'=e_j'(\mu)+1$ for brevity. We know that $w$ is either $w_j=e_{\nu}(x_j^{f_j}\otimes 1_{\mu})$ or $w_j'=e_{\nu'}(\eta_j^{f_j'}\otimes 1_{\mu})$ for some $j\in S(\mu)$ where $\nu=\mu_{f_ju_j,\uo}$ and $\nu'=\mu_{\uo,f'_ju_j}$ are weights of $w_j$ and $w_j'$, respectively. We consider two cases.

(a) Suppose $w=w_j$. Since $\mu=\lambda_{\us,\ut}$ we have 
$$\nu=\lambda_{\us,\ut}\widehat{a^{-f_ju_j}}\chi^{f_ju_j}=\lambda_{\us+f_ju_j,\ut}.$$ 
Therefore $\phi(w)=e_{\nu}x_j^{f_j}v_{\us,\ut}$ which by Lemma \ref{x-action} is seen to be 
$$e_{\nu}(\sum_{l+k=f_j}c_{l,k}v_{\us+lu_j,\ut+ku_j})\;\text{for some }\;c_{l,k}\in\k.$$
Furthermore, the equality of characters $\lambda_{\us+lu_j,\ut+ku_j}=\nu$ is equivalent to $a_j^{-l}b_j^k=a_j^{-f_j}$ which, in view of $f_j=l+k$, reduces to $(a_jb_j)^k=1$. By our assumption on datum the last condition holds only for $k=0$. Thus $e_{\nu}x_j^{f_j}v_{\us,\ut}=cv_{\us+f_ju_j,\ut},\;c\in\k$, and $c\neq 0$ by another application of Lemma \ref{x-action}. The last formula and part (i) make it clear that if $j\in S^{(1)}(\mu)=S^{(1)}(\lambda)$ and $s_j=0$, then $\phi(w)$ is  a primitive vector of rank $m+1$. If $s_j=e_j(\lambda)+ 1$, then $f_j=e_j({\mu})+1=m_j-e_j(\lambda)-1$ again by part (i), hence $f_j+s_j=m_j$, whence $\phi(w)=0$.

(b) We let $w=w_j'$. Proposition \ref{multipH*'}(1) and Lemma \ref{gconjug} show that $\eta_j^{f_j'}\otimes 1_{\mu}$ has weight $\nu'$, hence $w_j'=\eta_j^{f_j'}\otimes 1_{\mu}$. By Proposition \ref{eta&xiaction}(1) $\phi(w)=cv_{\us,\ut+f_j'u_j},c\neq 0$. If $j\in S^{(2)}(\lambda)$ and $t_j=0$, then $f_j'=e_j'(\lambda)+1$, hence $\phi(w)$ is a primitive vector of rank $m+1$. Otherwise, by an argument as above, $f_j'+t_j=m_j$, which gives $\phi(w)=0$. It remains to note that the primitive vectors constructed in (a) and (b) account for all primitive vectors of rank $m+1$.

For (4) we note that every element of $S\setminus S^{(3)}$ contributes a one and every element of $S^{(3)}$ contributes a two to the rank function. Thus the largest value the rank function can have is $|S\setminus S^{(3)}|+2|S^{(2)}|=|S^{(1)}|+|S^{(2)}|$.

(2) We begin the argument as in the same part of Theorem \ref{Lseries1}. We assume by the reverse induction on $m$ that $R^{m+1}=\Sigma_{\ell-m-1}$ and we let $M=I(\lambda)/R^{m+1}$. We want to prove the  equality $R^m/R^{m+1}=\Sigma(M)$. Assuming it does not hold, there is a simple $D$-module $L$  in $\Sigma(M)$ not contained in $R^m/R^{m+1}$. Let $k$ be the largest integer such that $L\subset R^k/R^{m+1}$. We denote by $\overline v$ the image of $v\in I(\lambda)$ in $M$. We define $\text{rk}(\overline{v_{\us,\ut}})=\text{rk}(v_{\us,\ut})$. 

Our first step is to show that $M$ is the free span of all $\overline{v_{\us,\ut}}$ with $\text{rk}(v_{\us,\ut})<m$. Indeed, a $v_{\us,\ut}$ of rank $\ge m+1$ is characterized by the property that $s_j\ge e_j+1$ or $t_{j'}\ge e_{j'}'+1$ for some $j,j'\in S(\lambda)$ of the total number greater than $m$. Therefore there is a primitive $v_{\uc,\ud}$ such that $s_j=c_j+k_j$ and $t_{j'}=d_{j'}+l_{j'}$ for some positive integers $k_j$ and $l_{j'}$. Now a glance at Propositions \ref{eta&xiaction} and \ref{x-action} leads to conclusion that $v_{\us,\ut}$ lies in $x^{\ui}\eta^{\uj}v_{\uc,\ud}$. Since by part (1) $R^{m+1}$ is generated by primitive vectors of rank $m+1$, the assertion follows.

Let $u$ be a generator of $L$ written in basis $B$ as
\begin{equation*}u=\sum c_{\us,\ut}\overline{v_{\us,\ut}},\,0\neq c_{\us,\ut}\in\k.\tag{*}\end{equation*}
The sum $u_k$ of all $\overline{v_{\uc,\ud}}$ of rank $k$ occuring in (*) is nonzero. Let's call the number of terms  in the sum (*) for $u_k$ {\em the length} of $u_k$. We pick a generator $u$ with $u_k$ of the smallest length. Fix one $\overline{v_{\uc,\ud}}$ involved in $u_k$. Assuming $k<m$ we can find $j\in S(\lambda)$ such that $c_j< e_j(\lambda)+1$ or $d_j< e_j'(\lambda)+1$. Suppose $c_j< e_j(\lambda)+1$. Set $l_j=e_j(\lambda)+1-c_j$ and $\nu=\lambda_{\uc+l_ju_j,\ud}$. By the argument used in the proof of (1,(iia)) for every vector $v_{\us,\ut}$, $e_{\nu}x_j^{l_j}v_{\us,\ut}=\kappa v_{\us+l_ju_j,\ut},\,\kappa\in\Bbbk^{\bullet}$. Since $v_{\uc+l_ju_j,\ut}$ has rank $k+1$, we see that $e_{\nu}x_j^{l_j}u$ is a nonzero generator with a lesser number of basis monomials of rank $k$, a contradiction.

Assuming $d_j< e_j'(\lambda)+1$, set ${l_j}'=e_j'(\lambda)+1- d_j$. Using the argument of part (1,(iib)) we deduce that $\eta^{{l_j'}'}u$ is a nonzero element of $L$ with the lesser number of basis monomials of rank $k$. This completes the proof of (2).

(3) We keep notation of part (2). By part (1) 
\begin{equation*}{\mathcal L}^m=\sum D.\overline{v_{\us,\ut}}\tag{**}\end{equation*} 
where $v_{\us,\ut}$ runs over all primitive basis vectors of rank $m$. Fix one $\overline{v_{\us,\ut}}$. By Propositions \ref{eta&xiaction} and \ref{x-action} $D.\overline{v_{\us,\ut}}$ is the span of the set\newline $B_{\us,\ut}=\{\overline{v_{\uc,\ud}}|\text{rk}(v_{\uc,\ud})=m\;\text{and}\;\uc\ge\us,\;\ud\ge\ut\}$. The condition $\text{rk}(v_{\uc,\ud})=m$ forces $0\le c_j\le e_j$ if $j\in S^{(1)}(\lambda)$ and $s_j=0$, and, likewise, $0\le d_j\le e'_j$ if $j\in S^{(2)}(\lambda)$ with $t_j=0$. If $j\in S^{(1)}(\lambda)$ or $j\in S^{(2)}(\lambda)$ with $s_j=e_j+1$ or $t_j=e'_j+1$, then $0\le c_j\le m_j-e_j-2$ or $0\le d_j\le m_j-e'_j-2$, respectively. For $j\notin S(\lambda),\,c_j,d_j$ take on every value in $[m_j]$. Therefore by part 1(i) and the dimension formula of Corollary \ref{dimension} $|B_{\us,\ut}|=\dim\,L(\lambda_{\us,\ut})$. As $D.\overline{v_{\us,\ut}}\supset L(\lambda_{\us,\ut})$ we obtain the equality $D.\overline{v_{\us,\ut}}=L(\lambda_{\us,\ut})$. Now, were sum (**) not direct, some $\overline{v_{\us,\ut}}$ would be a linear combimation of elements of other $B_{\us',\ut'}$. Since $\overline{v_{\us,\ut}}$ is a basis element, it would lie in some $D.\overline{v_{\us',\ut'}}$. However, from Theorem \ref{simpleD-mod} one sees that every simple $D$-module has a unique line of primitive elements, a contradiction.\qed

%Next we claim that every primitive vector of $M$ is a linear combination of the images of primitive vectors. Let $B$ be the basis of $M$ made up by the images of $v_{\us,\ut}$. Pick primitive vector $u$ of $M$ and expand $u$ relative to $B$,
%\begin{equation*}u=\sum c_{\us,\ut}\overline{v_{\us,\ut}},\,0\neq c_{\us,\ut}\in\k.\tag{*}\end{equation*}
%Proposition \ref{eta&xiaction}(2) shows that for every $k$ $\xi_k$ maps elements of $B$ into scalar multiples of elements of $B$. Therefore $\xi_ku=0$ is equivalent to $\xi_k\overline{v_{\us,\ut}}=0$ for all $(\us,\ut)$ occuring in (*). Since $\overline{v_{\us-u_k,\ut}}\in B$, same proposition gives that $\xi_k\overline{v_{\us,\ut}}=0$ if either $s_k=0$ or $\lambda(a_k\chi_k^{-1})=q_k^{-(s_k-1)}$. Therefore we deduce from Proposition \ref{y-action} that multiplication by $y_k$ also sends elements of $B$ into scalar multiples of elements of $B$. Thus $y_k\overline{v_{\us,\ut}}=0$ holds for all $(\us,\ut)$ and all $k$. Moreover, Proposition \ref{y-action} shows that $y_k\overline{v_{\us,\ut}}=0$ implies that either $t_k=0$ or $\lambda(b_k\chi_k)=q_k^{-(t_k-1)}$. Thus each $v_{\us,\ut}$ is primitive.

The (neo)classical quantum groups of Drinfel'd, Jimbo and Lusztig have the group of grouplike equal to the direct sum of cyclic subgroups generated by the grouplike associated to the `positive' or `negative' half of skew primitive generators. In the finite-dimensional case (see e.g. \cite{Lu2}) the orders of all $q_i$ are odd. In general those two conditions on datum are independent of each other. We call datum $\mathcal D$ {\em classical} if either all $|q_i|$ are odd, or the elements $\{a_i,b_i\}_{i\in\un}$ are independent in the sense that they generate subgroup equal to the direct sum of cyclic subgroups generated by $a_i$ and $b_i$. We turn to liftings $H$ with classical data. We will give a complete description of the lattice of submodules of $I(\lambda)$ for every $\lambda\in\GG$. This is possible because the lattice of $D$-submodules turns out to be distributive, a consequence of the next lemma.
\begin{Lem} Suppose $\mathcal D$ is classical. Then $\mathcal D$ is half-clean and the weights $\lambda_{\us,\ut},\us,\ut\in\mathbb Z(\um)$ are distinct.
\end{Lem}
\pf The first assertion holds by definition if all $|q_i|$ are odd. Else, suppose the set $\{a_i,b_i\}_{i\in\un}$ is independent. Then $\prod_{i=1}^n(a_ib_i)^t_i=1$ implies $a_i^{t_i}=1$ for all $i$. Since $\chi_i(a_i)=q_i$ and $q_i$ has order $m_i$, $s$ is divisible by $m_i$, which proves that $\mathcal D$ is half-clean.

Moving on to the second claim we must show that $\us=\us'$ and $\ut=\ut'$ whenever $\widehat{a^{-\us}b^{\ut}}\chi^{\us+\ut}=\widehat{a^{-\us'}b^{\ut'}}\chi^{\us'+\ut'}$. This equation is equivalent to
\begin{align*}a^{-\us}b^{\ut}=a^{-\us'}b^{\ut'}\\
\chi^{\us+\ut}=\chi^{\us'+\ut'}\end{align*}
Set $p_i=s_i-s'_i$ and $r_i=t_i-t_i'$ for all $i$. We rewrite the above two equations as 
\begin{align*}a^{-\up}b^{\ur}=1\tag{*}\\
\chi^{\up+\ur}=1\tag{**}\end{align*}        
Suppose $\{a_i,b_i\}_{i\in\un}$ are independent. The equation (*) implies equalities $a_i^{p_i}=1$ and $b_i^{r_i}=1$ for all $i$. As $\chi_i(a_i)=\chi_i(b_i)=q_i$ and the latter has order $m_i$, we see that $p_i$ and $r_i$ are divisible by $m_i$. Since $-m_i< p_i,r_i<m_i$ we conclude that $p_i=0=r_i$, and this holds for all $i$.

Next assume that all $m_i$ are odd. We induct on $n$ assuming the lemma holds for every datum on $<n$ points. Since for $\un=\emptyset$ the claim is vacuously true we proceed to the induction step. Applying $\chi_1$ to the equality (*) gives
\begin{equation*}\chi_1(a_1)^{-p_1}\chi_1(b_1)^{r_1}\prod_{i=2}^n\,\chi_1(a_i)^{-p_i}\chi_1(b_i)^{r_i}=1\tag{!}\end{equation*}
Using datum conditions (D1)-(D2) we have $\chi_1(a_1)=\chi_1(b_1)=q_1,\,\chi_1(a_i)^{-1}=\chi_i(a_1)$ and $\chi_1(b_i)=\chi_i(a_1)$. Therefore equality (!) takes on the form 
\begin{equation*}q_1^{r_1-p_1}\prod_{i=2}^n\chi_i(a_1)^{p_i+r_i}=1.\end{equation*}
Further, evaluating the left side of (**) at $a_1$ we get the equality
\begin{equation*}q_1^{p_1+r_1}\prod_{i=2}^n\chi_i(a_1)^{p_i+r_i}=1.\end{equation*}
It follows that $r_1-p_1\equiv r_1+p_1 \mod m_i$. In addition taking the value of the left side of (**) at $b_1$ we have
$$q_1^{p_1+r_1}\prod_{i=2}^n\chi_i(b_1)^{p_i+r_i}=q_1^{p_1+r_1}\prod_{i=2}^n\chi_i(a_1)^{-(p_i+r_i)}=1.$$
Comparing the last two equalities we see that $q_1^{2(p_1+r_1)}=1$ whence $p_1+r_1\equiv 0 \mod m_i$, as $m_i$ is odd. It follows that both $p_i$ and $r_i$ are divisible by $m_i$, hence $p_1=0=r_1$, and the proof is complete.\qed

The above lemma makes it clear that all vectors $v_{\us,\ut},\us,\ut\in\mathbb Z(\um)$ have distinct weights. Therefore every weight subspace $I_{\mu}(\lambda):=e_{\mu}I(\lambda)$ is one-dimensional. Since $I(\lambda)$ is a semisimple $\k\GG$-module $I(\lambda)\simeq \oplus_{\mu\in\GG}m_{\mu}L(\mu)$ where $m_{\mu}$ is the multiplicity of $L(\mu)$ in $I(\lambda)$. As $L(\mu)$ contains a $\mu$-weight vector, $m_{\mu}\le 1$ for all $\mu$. Thus $I(\lambda)$ is a multiplicity free module for all $\lambda\in\GG$.
We digress briefly into a general theory of such modules. For an alternate treatment see \cite{Al}

Let $A$ be an algebra and $M$ a left $A$-module of finite length with every simple $A$-module occuring at most once in a composition series of $M$. Let $\Lambda$ be the submodule lattice of $M$. $\Lambda$ is distributive by a standard criterion \cite[II.13]{B}. An element $J\neq 0$ will be called {\em local} (or join-irreducible \cite{B}) if $A\underset\neq\subset J$ and $B\underset\neq\subset J$ imply $A+B\underset\neq\subset J$. Clearly the radical $R(J)$ is a unique maximal submodule of $J$. We let $R(J)=0$ if $J$ is simple. Let $\mathcal J={\mathcal J}(M)$ denote the poset (partially ordered set) of local submodules ordered by inclusion. By \cite[Cor. III.3]{B} $\mathcal J$ forms a distinguished basis for $\Lambda$ in the sense that every $X\in\Lambda$ has a unique representation as the sum of an irredundant set of local submodules. Let $\mathcal S$ be the set of all composition factors of $M$. It turns out that $\mathcal J$ also determines $\mathcal S$ and the way composition factors are `stuck' together. The precise statement is-
\begin{Prop}\label{distrlattice} In the foregoing notation
\begin{itemize}\item[(1)] Let $N=J_1+\cdots +J_k$ be an irredundant sum of local submodules of $M$. Then each $J_i$ is a maximal submodule of $N$, $$R(N)=\sum_{i=1}^kR(J_i)\;\text{and}\;N/R(N)=\bigoplus_{i=1}^kJ_i/R(J_i).$$
\item[(2)] The mapping $J\mapsto J/R(J)$ sets up a bijection between $\mathcal J$ and~ $\mathcal S$.
\item[(3)] The composition length of $M$ equals $|\mathcal J|$.
\item[(4)] For $L\in\mathcal S$ let $J(L)$ be the preimage of $L$ under the map in {\rm(2)}. For any two simple modules $L$ and $L'$, $L$ occurs before $L'$ in a composition series of $M$ if and only if $J(L)$ and $J(L')$ are either incomparable in $\mathcal J$ or $J(L)\supset J(L')$.
\end{itemize}
\end{Prop}
\pf (1) Suppose $N$ is as in (1). Let $K=\sum_{i=1}^kR(J_i)$. Then $J_i\cap K=R(J_i)+\sum_{k\neq i}J_i\cap R(J_k)=R(J_i)$ because $J_i\cap R(J_k)\subset J_i\cap J_k\subset R(J_i)$. Therefore $N/K\simeq \bigoplus_{i=1}^k J_i/R(J_i)$, hence $K\supset R(N)$. On the other hand, ${J_i+R(N)}/R(N)$ is semisimple for all $i$, hence $R(N)\cap J_i\subset R(J_i)$, whence $R(N)\subset K$.

(2) Taking filtration $M\supset R(M)\supset R^2(M)\supset\cdots \supset 0$ we see by part (1) that every composition factor of $M$ is of the form $J/R(J)$ for some $J\in\mathcal J$. Thus the mapping $\mathcal J\to\mathcal S,\;J\mapsto J/R(J)$ is onto. Were $J/R(J)\simeq J'/R(J')$ for some $J\neq J'$, the multiplicity $m$ of $J/R(J)$ in $M$ would be $\ge 2$. For, if $J\supset J'$, then refining $J\supset J'\supset 0$ we get $m\ge 2$. Else, we refine $J+J'\supset J\supset 0$, and get $m\ge 2$, again. By our assumption that $M$ is multiplicity free, the assertion follows. 

(3) is \cite[Lemma III.2]{B}. It also follows immediately from part (2) as the composition length is $|\mathcal S|=|\mathcal J|$.

(4) Suppose there is a composition series of $M$ with $L$ preceding $L'$. Say $L=A/B$ and $L'=B/C$. By part (1) $J(L)\subset A$ and $J(L')\subset B$. Since $J(L)\not\subset B$, lest we have the multiplicity of $L$ $\ge 2$, we see that $J(L)\not\subset J(L')$. Suppose $J(L)$ and $J(L')$ are incomparable. The refining $M\supset J(L)+J(L')\supset J(L')\supset 0$ we get $L$ before $L'$, and refining $M\supset J(L)+J(l')\supset J(L)\supset 0$ reverses the order of their appearance.\qed

We return to modules $I(\lambda)$ assuming the datum to be classical. For every primitive vector $v_{\us,\ut}$ we define pair of sets $(S^{\prime}(v_{\us,\ut}),S^{\prime\prime}(v_{\us,\ut})$ by $S^{\prime}(v_{\us,\ut})=\{j\in S({\lambda})|s_j=e_j+1\}$ and $S^{\prime\prime}(v_{\us,\ut})=\{j\in S({\lambda})|t_j=e_j'+1\}$. We let $P=P(\lambda)$ denote the set of all pairs $(S^{\prime}, S^{\prime\prime})$ with $S^{\prime}\subset S^{\prime}({\lambda})$ and $S^{\prime\prime}\subset S^{\prime\prime}(\lambda)$. We turn $P$ into a poset by defining an ordering $(S^{\prime},S^{\prime\prime})\succeq (T^{\prime}, T^{\prime\prime})$ if and only if $S^{\prime}\subseteq T^{\prime}$ and $S^{\prime\prime}\subseteq T^{\prime\prime}$. The main result is 
\begin{Thm}\label{classicalliftings} In the foregoing notations with $H=H(\mathcal D)$ and $D=D(H)$
\begin{itemize}\item[(1)] A submodule $J$ of $I(\lambda)$ is local if and only if $J$ is generated by a primitive weight vector.
\item[(2)] The poset $\mathcal J$ of local submodules of $I(\lambda)$ is isomorphic to $P$.
\item[(3)] The composition length of $I(\lambda)$ equals $2^{\ell(\lambda)-1}$
\end{itemize}
\end{Thm}
\pf (1) Since all $v_{\us,\ut}$ have distinct weights every submodule $M$ of $I(\lambda)$ is the span of $M\cap B$ where $B=\{v_{\us,\ut}|\us,\ut\in\mathbb Z(\um)\}$. Therefore, if $J$ is generated by a single basis vector $v\in B$, then were $J=K+M$, $v$ would lie in $K$ or $N$, hence $J$ is local. Conversely, if $J$ is local, then, as $J=\sum_{v\in J\cap B}Dv$, we have $J=Dv$ for some $v$.

Let $v=v_{\us,\ut}$ be a generator of $J$. If $s_k\neq 0$ for some $k\notin S(\lambda)$ or $s_k\neq e_k+1$ for some $k\in S(\lambda)$, then by Propositions \ref{eta&xiaction} and \ref{x-action} $\xi_kv_{\us,\ut}$ also generates $J$. Thus we can assume that $s_k=0$ for every $k\notin S(\lambda)$ and $s_k=0,e_k+1$ for every $k\in S(\lambda)$. We come to a similar conclusion about every $t_k$, viz. $t_k=0$ if $k\notin S(\lambda)$ and $t_k=0,e_k'+1$ for every $k\in S(\lambda)$ by using Propositions \ref{eta&xiaction} and \ref{y-action}. This proves (1)

(2) Let $\psi:\mathcal J\to P$ be the mapping sending $J$ to $\psi(v):=(S^{\prime}(v),S^{\prime\prime}(v))$ where $v$ is a primitive generator of $J$. By  primitivity of $v$, $J$ is the span of all $x^{\ui}\eta^{\uj}v$. Therefore Propositions \ref{eta&xiaction}(2) and \ref{x-action} give the relation $\psi(w)\preceq\psi(v)$ for every primitive weight vector $w\in J$. This shows first that if $v$ and $w$ generate $J$ then $\psi(w)=\psi(v)$, so that
$\psi$ is well-defined, and second that $\psi$ is isotonic, which yields (2)

(3) By part (3) of the preceding proposition the composition length of $I(\lambda)$ equals $|\mathcal J|$, hence $|P|$. The latter is $2^{|S^{(1)}|+|S^{(2)}|}$ which yields the claim by Theorem \ref{Lseries2}(4).\qed

\end{document}